\documentclass[10pt,onecolumn]{amsart}

\usepackage{url}
\usepackage[usenames, dvipsnames]{color}
\usepackage[colorlinks=true,
        raiselinks=true,
        linkcolor=MidnightBlue,
        citecolor=Mahogany,
        urlcolor=ForestGreen,
        pdfauthor=Kostas Margellos,
        pdftitle={Hamilton-Jacobi formulation for Reach-Avoid Differential Games},
        pdfkeywords={Reachability, differential game theory, optimal control, Hamilton-Jacobi equations, hybrid systems},
        pdfsubject={Technical Report},
        plainpages=false]{hyperref}
\usepackage{amsmath, graphicx, fancyhdr, graphics}
\usepackage{subfigure}
\usepackage{algorithm}
\usepackage{algorithmic}

\title{Hamilton-Jacobi formulation for Reach-Avoid Differential Games}
\author[K.~Margellos]{Kostas~Margellos}
\author[J.~Lygeros]{John~Lygeros}
\thanks{K. Margellos and J. Lygeros are with the Automatic Control Laboratory, Department
of Electrical Engineering and Information Technology, Swiss Federal Institute of Technology (ETH),
Physikstrasse 3, 8092 Z\"urich, e-mail: {margellos, lygeros}@control.ee.ethz.ch}

\begin{document}
\maketitle
\begin{abstract}
A new framework for formulating reachability
problems with competing inputs, nonlinear dynamics and state constraints as optimal
control problems is developed. Such reach-avoid problems
arise in, among others, the study of safety problems in hybrid
systems. Earlier approaches to reach-avoid computations are
either restricted to linear systems, or face numerical difficulties
due to possible discontinuities in the Hamiltonian of the optimal
control problem. The main advantage of the approach proposed
in this paper is that it can be applied to a general class of target hitting continuous dynamic games with nonlinear dynamics,
and has very good properties in terms of its numerical solution,
since the value function and the Hamiltonian of the system are
both continuous. The performance of the proposed method is
demonstrated by applying it to a two aircraft collision avoidance scenario under target window constraints and in the presence of wind disturbance. Target Windows are a novel concept in air traffic management, and represent spatial and temporal constraints, that the aircraft have to respect to meet their
schedule.
\end{abstract}

\section{Introduction}
Reachability for continuous and hybrid systems has been an important topic of research in the
dynamics and control literature. Numerous problems regarding safety of air traffic
management systems \cite{Tomlin_Lygeros1}, \cite{Tomlin_Pappas}, flight control
\cite{Tomlin_Lygeros2}, \cite{Tomlin_Lygeros3}, \cite{Kitsios} ground transportation systems \cite{Livadas}, \cite{Lygeros_Godbole}, etc. have been formulated in the framework of reachability theory.
In most of these applications the main aim was to design suitable controllers to steer or keep the state of the system in a "safe" part of the state space. The synthesis of such safe controllers for hybrid systems relies on the ability to solve target problems for the case where state constraints are also present. The sets that represent the solution to those problems are known as capture basins \cite{Aubin}.
One direct way of computing these sets was proposed in \cite{Gardaliaguet1}, \cite{Gardaliaguet2}, and was formulated in the context of viability theory \cite{Aubin}. Following the same approach, the authors of \cite{ALQSS}, \cite{Gao_Lygeros} formulated viability, invariance and pursuit-evasion gaming problems for hybrid systems and used non-smooth analysis tools to characterize their solutions.
Computational tools to support this approach have been already developed by \cite{Saint-Pierre}.

An alternative, indirect way of characterizing such problems is through the level sets of the value function of an appropriate optimal control problem. By using dynamic programming, for reachability/invariant/viability problems without state constraints, the value function can be characterized as the viscosity solution to a first order
partial differential equation in the standard Hamilton-Jacobi form \cite{Lygeros},
\cite{Evans_Souganidis}, and \cite{Mitchell_Bayen1}. Numerical algorithms based on level set methods have been developed by \cite{Osher_Sethian}, \cite{Sethian}, have been coded in efficient computational tools by \cite{Mitchell_Bayen1}, \cite{Mitchell_Tomlin} and can be directly applied to reachability computations.

In the case where state constraints are also present, this target hitting problem is the solution to a reach-avoid problem in the sense of \cite{Tomlin_Lygeros1}. The authors of \cite{Tomlin_Lygeros1}, \cite{TomlinPhD} developed a reach-avoid computation, whose value function was characterized as a solution to a pair of coupled variational inequalities. In \cite{Mitchell_Tomlin}, \cite{MitchellPhD}, \cite{Oishi} the authors proposed another characterization, which involved only one Hamilton-Jacobi type partial differential equation together with an inequality constraint. These methods are hampered from a numerical computation point of view by the fact that the Hamiltonian of the system is in general discontinuous \cite{TomlinPhD}.

In \cite{Kurzhanski_Varaiya1}, a scheme based on ellipsoidal techniques so as to compute reachable sets for control systems with constraints on the state was proposed. This approach was restricted to the class of linear systems. In \cite{Kurzhanski_Varaiya2}, this approach was extended to a list of interesting target problems with state constraints. The calculation of a solution to the equations proposed in \cite{Kurzhanski_Varaiya1}, \cite{Kurzhanski_Varaiya2} is in general not easy apart from the case of linear systems, where duality techniques of convex analysis can be used.

In this paper we propose a new framework of characterizing reach-avoid sets of nonlinear control systems as the solution to an optimal control problem. We consider the case where we have competing inputs and hence adopt the gaming formulation proposed in \cite{Evans_Souganidis}. We first restrict our attention to a specific reach-avoid scenario, where the objective of the control input is to make the states of the system hit the target at the end of our time horizon and without violating the state constraints, while the disturbance input tries to steer the trajectories of the system away from the target. We then generalize our approach to the case where the controller aims to steer the system towards the target not necessarily at the terminal, but at some time within the specified time horizon. Both problems could be treated as pursuit-evasion games, and for a worst case setting we define a value function similar to \cite{Kurzhanski_Varaiya2} and prove that it is the unique continuous viscosity solution to a quasi-variational inequality of a form similar to \cite{Georgiou}, \cite{Barron_Ishii}. The advantage of this approach is that the properties of the value function and the Hamiltonian (both of them are continuous) enable us using existing tools to compute the solution of the problem numerically.

To illustrate our approach, we consider a reach-avoid problem that arises in the area of air traffic management, in particular the problem of collision avoidance in the presence of 4D constraints, called Target Windows. Target Windows (TW) are spatial and temporal constraints and form the basis of the CATS research project \cite{CATS}, whose aim is to increase punctuality and predictability during the flight. In \cite{CDCpaper} a reachability approach of encoding TW constraints was proposed. We adopt this framework and consider a multi-agent setting, where each aircraft should respect its TW constraints while avoiding conflict with other aircraft in the presence of wind. Since both control and disturbance inputs (in our case the wind) are present, this problem can be treated as a pursuit-evasion differential game with state constraints, which are determined dynamically by performing conflict detection.

In Section II we pose two reach-avoid problems for continuous systems with competing inputs and state constraints, and formulate them in the optimal control framework. Section III provides the characterization of the value functions of these problems as the viscosity solution to two variational inequalities. In Section IV we present an application of this approach to a two aircraft collision avoidance scenario with realistic data. Finally, in Section V we provide some concluding remarks and directions for future work.
\section{Differential games and Reach-Avoid problems}

\subsection{Differential game problem formulation}
Consider the continuous time control system $\dot{x}=f(x,u,v)$, and an arbitrary time horizon $T \geq 0$.
with $x \in \mathbb{R}^n$, $u \in \textit{U} \subseteq \mathbb{R}^m$, $v \in \textit{V} \subseteq \mathbb{R}^p$, and $f(\cdot,\cdot, \cdot): \mathbb{R}^n \times U \times V \rightarrow \mathbb{R}^n$.
Let $\mathcal{U}_{[t,t']}$, $\mathcal{V}_{[t,t']}$ denote the set of Lebesgue measurable functions from the interval $[t,t']$ to $\textit{U}$, and $\textit{V}$ respectively. Consider also two functions $l(\cdot): \mathbb{R}^n \rightarrow \mathbb{R}$,  $h(\cdot): \mathbb{R}^n \rightarrow \mathbb{R}$ to be used to encode the target and state constraints respectively,\\
\newline
\textbf{Assumption 1.} \textit{ $\textit{U} \subseteq \mathbb{R}^m$ and $\textit{V} \subseteq \mathbb{R}^p$ are compact. $f$, $l$ and $h$ are bounded and Lipschitz continuous in x and continuous in u and v.}

Under Assumption $1$ the system admits a unique solution $x(\cdot): [t,T] \rightarrow \mathbb{R}^n$  for all $t \in [0,T]$, $u(\cdot) \in \mathcal{U}_{[t,T]}$ and $v(\cdot) \in \mathcal{V}_{[t,T]}$. For $\tau \in [t,T]$ this solution will be denoted as
\begin{equation}\label{phi}
    \phi(\tau,t,x,u(\cdot),v(\cdot))=x(\tau).
\end{equation}
Let $C_f>0$ be a bound such that for all $x,\hat{x} \in \mathbb{R}^n$ and $u(\cdot) \in \mathcal{U}_{[t,T]}$ and for all $u \in U$,
\begin{equation*}
|f(x,u)| \leq C_f \text{ and } |f(x,u)-f(\hat{x},u)| \leq C_f|x-\hat{x}|.
\end{equation*}
Let also $C_l>0$ and $C_h>0$ be such that
\begin{align*}
|l(x)| \leq C_l &\text{ and } |l(x)-l(\hat{x})| \leq C_l|x-\hat{x}|, \\
|h(x)| \leq C_h &\text{ and } |h(x)-h(\hat{x})| \leq C_h|x-\hat{x}|.
\end{align*}
In a game setting it is essential to define the information patterns that the two players use. Following \cite{P.P.Varaiya}, \cite{Evans_Souganidis} we restrict the first player to play non-anticipative strategies. A non-anticipative strategy is a function $\gamma: \mathcal{V}_{[0,T]} \rightarrow \mathcal{U}_{[0,T]}$ such that for all $s \in [t,T]$ and for all $v,\hat{v} \in \mathcal{V}$, if $v(\tau)=\hat{v}(\tau)$ for almost every $\tau \in [t,s]$, then $\gamma[v](\tau)=\gamma[\hat{v}](\tau)$ for almost every $\tau \in [t,s]$. We then use $\Gamma_{[t,T]}$ to denote the class of non-anticipative strategies.

Consider the sets $R$, $A$ related to the level sets of the two bounded, Lipschitz continuous functions $l: \mathbb{R}^n \rightarrow \mathbb{R}$ and  $h: \mathbb{R}^n \rightarrow \mathbb{R}$ respectively. For technical purposes assume that $R$ is closed whereas $A$ is open. Then $R$ and $A$ could be characterized as
\begin{equation*}
    R = \{x \in \mathbb{R}^n ~|~ l(x)\leq0\},~
    A = \{x \in \mathbb{R}^n ~|~ h(x) > 0\} .
\end{equation*}

\subsection{Reach-Avoid at the terminal time}
Consider now a closed set $R \subseteq \mathbb{R}^n$ that we would like to reach while avoiding an open set $A \subseteq \mathbb{R}^n$. One would like to characterize the set of the initial states from which trajectories can start and reach the set $R$ at the terminal time $T$ without passing through the set $A$ over the time horizon $[t,T]$. To answer this question on needs to determine whether there exists a choice of $\gamma \in \Gamma_{[t,T]}$ such that for all $v(\cdot) \in \mathcal{V}_{[t,T]}$, the trajectory $x(\cdot)$ satisfies $x(T) \in R$ and $x(\tau) \in A^c$ for all $\tau \in [t,T]$.

The set of initial conditions that have this property is then
\begin{align}
  RA(t,R,A) &= \{x \in \mathbb{R}^n ~|~ \exists \gamma(\cdot) \in \Gamma_{[t,T]}, ~\forall v(\cdot) \in \mathcal{V}_{[t,T]}, \\
  &( \phi(T,t,x,\gamma(\cdot),v(\cdot)) \in R) \land (\forall \tau \in [t,T], ~ \phi(\tau,t,x,\gamma(\cdot),v(\cdot)) \notin A ) \}. \nonumber
\end{align}

Now introduce the value function $V : \mathbb{R}^n \times [0,T] \rightarrow \mathbb{R}$
\begin{equation}
  V(x,t)=\inf_{\gamma(\cdot) \in  \Gamma_{[t,T]}}  \sup_{v(\cdot) \in  \mathcal{V}_{[t,T]}} \max  \{l(\phi(T,t,x,u(\cdot),v(\cdot))) , \max_{\tau \in [t,T]} h(\phi(\tau,t,x,u(\cdot),v(\cdot))) \}.
\end{equation}
$V$ can be thought of as the value function of a differential game, where $u$ is trying to minimize, whereas $v$ is trying to maximize the maximum between the value attained by $l$ at the end $T$ of the time horizon and the maximum value attained by $h$ along the state trajectory over the horizon $[t,T]$. Based on \cite{Lygeros}, \cite{Evans_Souganidis} and \cite{Georgiou}, we will show that the value function defined by $(3)$ is the unique viscosity solution of the following quasi-variational inequality.
\begin{equation}
    \max \{ h(x)-V(x,t), \frac{\partial V}{\partial t}(x,t)+\sup_{v \in \textit{V}} \inf_{u \in \textit{U}} \frac{\partial V}{\partial x}(x,t) f(x,u,v) \}=0,
\end{equation}
with terminal condition $V(x,T)= \max \{l(x),h(x) \}$.

It is then easy to link the set $RA(t,R,A)$ of $(2)$ to the level set of the value function $V(x,t)$ defined in $(3)$. \\
\newline
\textbf{Proposition 1.} $RA(t,R,A) = \{x \in \mathbb{R}^n ~|~ V(x,t)\leq0\}$. \\
\begin{proof}
$V(x,t)\leq0$ if and only if $\inf_{\gamma(\cdot) \in  \Gamma_{[t,T]}} \sup_{v(\cdot) \in  \mathcal{V}_{[t,T]}} \max \{l(\phi(T,t,x,\gamma(\cdot),v(\cdot))), \\ \max_{\tau \in [t,T]} h(\phi(\tau,t,x,\gamma(\cdot),v(\cdot) )) \} \leq 0$. Equivalently, there exists a strategy $\gamma(\cdot) \in \Gamma_{[t,T]}$ such that for all $v(\cdot) \in  \mathcal{V}_{[t,T]}$, $\max \{l(\phi(T,t,x,\gamma(\cdot),v(\cdot))) , \max_{\tau \in [t,T]} h(\phi(\tau,t,x,\gamma(\cdot),v(\cdot))) \} \leq 0$. The last statement is equivalent to there exists a $\gamma(\cdot) \in \Gamma_{[t,T]}$ such that for all $v(\cdot) \in  \mathcal{V}_{[t,T]}$, $l(\phi(T,t,x,\gamma(\cdot),v(\cdot))) \leq 0$ and $\max_{\tau \in [t,T]} h(\phi(\tau,t,x,\gamma(\cdot),v(\cdot))) \leq 0$. Or in other words, there exists a $\gamma(\cdot) \in \Gamma_{[t,T]}$ such that for all $v(\cdot) \in  \mathcal{V}_{[t,T]}$, $\phi(T,t,x,\gamma(\cdot),v(\cdot)) \in R$ and for all $\tau \in [t,T]~ \phi(\tau,t,x,\gamma(\cdot),v(\cdot)) \notin A$.
\end{proof}

\subsection{Reach-Avoid at any time}
Another related problem that one might need to characterize is the set of initial states from which trajectories can start, and for any disturbance input can reach the set $R$ not at the terminal, but at some time within the time horizon $[t,T]$, and without passing through the set $A$ until they hit $R$. In other words, we would like to determine the set
\begin{align}
 \widetilde{RA}(t,&R,A) = \{x \in \mathbb{R}^n ~|~ \exists \gamma(\cdot) \in \Gamma_{[t,T]}, ~\forall v(\cdot) \in \mathcal{V}_{[t,T]}, \\ &\exists \tau_1 \in [t,T],~
  ( \phi(\tau_1,t,x,\gamma(\cdot),v(\cdot)) \in R )\land (\forall \tau_2 \in [t,\tau_1], ~ \phi(\tau_2,t,x,\gamma(\cdot),v(\cdot)) \notin A ) \}. \nonumber
\end{align}

Based on \cite{Mitchell_Bayen2}, define the augmented input as $\tilde{u} = [u ~ \bar{u}] \in \textit{U} \times [0,1]$ and consider the dynamics
\begin{equation}
\tilde{f}(x,\tilde{u},v) = \bar{u}f(x,u,v).
\end{equation}
Let $\tilde{\phi}(\tau,x,t,\tilde{u}(\cdot),v(\cdot))$ denote the solution of the augmented system, and define $\widetilde{U}$, $\mathcal{\widetilde{U}}$ and $\widetilde{\Gamma}$ similarly to the previous case.
Following \cite{Mitchell_Bayen2} for every $\tilde{u} \in \mathcal{\widetilde{U}}_{[t,T]}$ the pseudo-time variable $\sigma : [t,T] \rightarrow [t,T]$ is given by
\begin{equation}
\sigma (\tau) = t + \int_{t}^{\tau} \bar{u}(s) ds.
\end{equation}
Consider $\sigma^*$ to be almost an inverse of $\sigma$ in the sense that $\sigma(\sigma^*(\tau))=\tau$. In \cite{Mitchell_Bayen2}, $\sigma^*$ was defined as the limit of a convergent sequence of functions, and it was shown that
\begin{equation}
\phi(\sigma(\tau), x,t,u(\sigma^*(\cdot)),v(\sigma^*(\cdot))) = \tilde{\phi}(\tau, x,t,\tilde{u}(\cdot),v(\cdot)) ,
\end{equation}
for any $\tau \in [t,T]$. Based on the analysis of \cite{Mitchell_Bayen2}, equation $(8)$ implies that the trajectory $\tilde{\phi}$ of the augmented system visits only the subset of the states visited by the trajectory $\phi$ of the original system in the time interval $[t,\sigma(\tau)]$.

Define now the value function
\begin{equation*}
  \widetilde{V}(x,t)=\inf_{\tilde{\gamma}(\cdot) \in  \widetilde{\Gamma}_{[t,T]}}  \sup_{v(\cdot) \in  \mathcal{V}_{[t,T]}} \max \{l(\tilde{\phi}(T,t,x,\tilde{\gamma}[v](\cdot),v(\cdot))) , \max_{\tau \in [t,T]} h(\tilde{\phi}(\tau,t,x,\tilde{\gamma}[v](\cdot),v(\cdot))) \}.
\end{equation*}
One can then show that $\widetilde{V}$ is related to the set $\widetilde{RA}$. \\
\newline
\textbf{Proposition 2.} \textit{ For $\tau \in [0,T]$, $\widetilde{RA}(\tau,R,A)= \{x \in \mathbb{R}^n ~|~ \widetilde{V}(x,\tau)\leq 0\}$. }

The proof of this proposition is given in Appendix A.

\section{Characterization of the value function}
\subsection{Basic properties of V}
We first establish the consequences of the principle of optimality for $V$. \\
\newline
\textbf{Lemma 1.} \textit{ For all $(x,t) \in \mathbb{R}^n \times [0,T]$ and all $\alpha \in [0,T-t]$:}
\begin{equation} V(x,t) = \inf_{\gamma(\cdot)\in  \Gamma_{[t,t+\alpha]}} \sup_{v(\cdot)\in  \mathcal{V}_{[t,t+\alpha]}} \big[ \max \big \{ \max_{\tau \in [t,t+\alpha]} h(\phi(\tau,t,x,u(\cdot))), V(\phi(t+\alpha,t,x,u(\cdot)),t+\alpha) \big \} \big]. \end{equation}
\textit{ Moreover, for all $(x,t) \in \mathbb{R}^n \times [0,T], ~ V(x,t) \geq h(x)$ .}

The proof for the second part is straightforward and follows from the definition of $V$. The proof for the first part is given in Appendix B.

We now show that $V$ is a bounded, Lipschitz continuous function.\\
\newline
\textbf{Lemma 2.} \textit{ There exists a constant $C > 0$ such that for all $(x,t),(\hat{x},\hat{t}) \in \mathbb{R}^n \times [0,T]$: }
\begin{equation*}
|V(x,t)| \leq C \text{ and } |V(x,t) - V(\hat{x},\hat{t})| \leq C(|x-\hat{x}|+|t-\hat{t}|).
\end{equation*}

The proof of this Lemma is given in Appendix B.

\subsection{Variational inequality for $V$}
We now introduce the Hamiltonian $H~:~ \mathbb{R}^n \times \mathbb{R}^n \rightarrow \mathbb{R}$, defined by
\begin{equation*}
H(p,x) = \sup_{v \in \textit{V}} \inf_{u \in \textit{U}}p^T f(x,u).
\end{equation*}
\newline
\textbf{Lemma 3.} \textit{There exists a constant $C > 0$ such that for all $p,q \in \mathbb{R}^n$, and all $x,y \in \mathbb{R}^n$}:
\begin{equation*}
|H(p,x)-H(q,x)| < C|p-q| \textit{ and } |H(p,x)-H(p,y)| < C|p||x-y|.
\end{equation*}

The proof of this fact is straightforward (see \cite{Lygeros}, \cite{Evans_Souganidis} for details). We are now in a position to state and prove the following Theorem, which is the main result of this section.\\
\newline
\textbf{Theorem 1.} \textit{ $V$ is the unique viscosity solution over $(x,t) \in \mathbb{R}^n \times [0,T]$ of the variational inequality}
\begin{equation*}
    \max \{ h(x)-V(x,t), \frac{\partial V}{\partial t}(x,t)+\sup_{v \in \textit{V}} \inf_{u \in \textit{U}} \frac{\partial V}{\partial x}(x,t) f(x,u,v) \}=0,
\end{equation*}
\textit{ with terminal condition $V(x,T)= \max \{l(x),h(x) \}$.}
\begin{proof}
Uniqueness follows from Lemma 2, Lemma 3 and \cite{Georgiou}. Note also that by definition of the value function we have $V(x,T)=\max \{l(x),h(x) \}$. Therefore it suffices to show that
\begin{enumerate}
\item For all $(x_0,t_0) \in \mathbb{R}^n \times (0,T)$ and for all smooth $W$ : $\mathbb{R}^n \rightarrow \mathbb{R}$, if $V-W$ attains a local maximum at $(x_0,t_0)$, then
   \begin{equation*}\max \{ h(x_0)-V(x_0,t_0), \frac{\partial W}{\partial t}(x_0,t_0)+\sup_{v \in \textit{V}} \inf_{u \in \textit{U}} \frac{\partial W}{\partial x}(x_0,t_0) f(x_0,u,v) \} \geq 0 \end{equation*}
\item For all $(x_0,t_0) \in \mathbb{R}^n \times (0,T)$ and for all smooth $W$ : $\mathbb{R}^n \rightarrow \mathbb{R}$, if $V-W$ attains a local minimum at $(x_0,t_0)$, then
    \begin{equation*}\max \{ h(x_0)-V(x_0,t_0), \frac{\partial W}{\partial t}(x_0,t_0)+\sup_{v \in \textit{V}} \inf_{u \in \textit{U}} \frac{\partial W}{\partial x}(x_0,t_0) f(x_0,u,v) \} \leq 0 \end{equation*}
\end{enumerate}
The case $t=0$ is automatically captured by \cite{Evans}.\\
\textbf{Part 1.}
Consider an arbitrary $(x_0,t_0) \in \mathbb{R}^n \times (0,T)$ and a smooth $W : \mathbb{R}^n \times (0,T) \rightarrow \mathbb{R}$ such that $V-W$ has a local maximum at $(x_0,t_0)$. Then, there exists $\delta_1 > 0$ such that for all $(x,t) \in \mathbb{R}^n \times (0,T)$ with $|x-x_0|^2 + (t-t_0)^2 < \delta_1$
\begin{equation*}
(V-W)(x_0,t_0) \geq (V-W)(x,t).
\end{equation*}
We would like to show that
\begin{equation*}
    \max \{ h(x_0)-V(x_0,t_0), \frac{\partial W}{\partial t}(x_0,t_0)+\sup_{v \in \textit{V}} \inf_{u \in \textit{U}} \frac{\partial W}{\partial x}(x_0,t_0) f(x_0,u,v) \} \geq 0.
\end{equation*}
Since by Lemma 1 $h(x)-V(x,t) \leq 0$, either $h(x_0)=V(x_0,t_0)$ or, $h(x_0)-V(x_0,t_0) < 0$. For the former the claim holds, whereas for the latter it suffices to show that there exists $v \in \textit{V}$ such that for all $u \in \textit{U}$
\begin{equation*}
\frac{\partial W}{\partial t}(x_0,t_0)+\frac{\partial W}{\partial x}(x_0,t_0)f(x_0,u,v) \geq 0.
\end{equation*}
For the sake of contradiction assume that for all $v \in \textit{V}$ there exists $u \in \textit{U}$ such that for some $\theta > 0$
\begin{equation*}
\frac{\partial W}{\partial t}(x_0,t_0)+\frac{\partial W}{\partial x}(x_0,t_0)f(x_0,u,v) < -2\theta < 0.
\end{equation*}
Since $W$ is smooth and $f$ is continuous, then based on \cite{Evans_Souganidis} we have that
\begin{equation*}
\frac{\partial W}{\partial t}(x_0,t_0)+\frac{\partial W}{\partial x}(x_0,t_0)f(x_0,u,\zeta) < -\frac{3\theta}{2} < 0,
\end{equation*}
for all $\zeta \in B(v,r) \cap \textit{V}$ and some $r>0$, where $B(v,r)$ denotes a ball centered at $v$ with radius $r$. Because $\textit{V}$ is compact there exist finitely many distinct points $v_1,...,v_n \in \textit{V}, ~~ u_1,...,u_n \in \textit{U}$, and $r_1,...,r_n >0$ such that $\textit{V} \subset \bigcup_{i=1}^{n} B(v_i,r_i)$ and for $\zeta \in B(v_i,r_i)$
\begin{equation*}
\frac{\partial W}{\partial t}(x_0,t_0)+\frac{\partial W}{\partial x}(x_0,t_0)f(x_0,u_i,\zeta) < -\frac{3\theta}{2} < 0.
\end{equation*}
Define $g:~\textit{V} \rightarrow \textit{U}$ by setting for $k=1,...,n$, $g(v)=u_k$ if $v \in B(u_k,r_k) \backslash \bigcup_{i=1}^{k-1} B(u_i,r_i)$. Then
\begin{equation*}
\frac{\partial W}{\partial t}(x_0,t_0)+\frac{\partial W}{\partial x}(x_0,t_0)f(x_0,g(v),v) < -\frac{3\theta}{2} < 0.
\end{equation*}
Since $W$ is smooth and $f$ is continuous, there exists $\delta_2 \in (0,\delta_1)$  such that for all $(x,t) \in \mathbb{R}^n \times (0,T)$ with $|x-x_0|^2 + (t-t_0)^2 < \delta_2$
\begin{equation*}
\frac{\partial W}{\partial t}(x,t)+\frac{\partial W}{\partial x}(x,t)f(x,g(v),v) < -\theta < 0.
\end{equation*}
Finally, define $\gamma:~\mathcal{V}_{[t_0,T]} \rightarrow \mathcal{U}_{[t_0,T]}$ by $\gamma[v](\tau) = g(v(\tau))$ for all $\tau \in [t_0,T]$. It is easy to see that $\gamma$ is now non-anticipative and hence $\gamma(\cdot) \in \Gamma_{[t_0,T]}$. So for all $v(\cdot) \in \textit{V}_{[t_0,T]}$ and all $(x,t) \in \mathbb{R}^n \times (0,T)$ such that $|x-x_0|^2 + (t-t_0)^2 < \delta_2$,
\begin{equation*}
\frac{\partial W}{\partial t}(x,t)+\frac{\partial W}{\partial x}(x,t)f(x,\gamma[v](\cdot),v(\cdot)) < -\theta < 0.
\end{equation*}
By continuity, there exists $\delta_3 > 0$  such that $|\phi(t,t_0,x_0,\gamma[v](\cdot),v(\cdot))-x_0|^2 + (t-t_0)^2 < \delta_2$ for all $t \in [t_0,t_0+\delta_3]$.
Therefore, for all $v(\cdot) \in \textit{V}_{[t_0,T]}$
\begin{align*}
V(\phi&(t, t_0,x_0,\gamma(\cdot),v(\cdot)),t)-V(x_0,t_0) \leq W(\phi(t, t_0,x_0,\gamma(\cdot),v(\cdot)),t)-W(x_0,t_0) \\
& = \int_{t_0}^{t} \Big( \frac{\partial W}{\partial s}(\phi(s, t_0,x_0,\gamma(\cdot),v(\cdot)),s)~ \\&+ \frac{\partial W}{\partial x}(\phi(s, t_0,x_0,\gamma(\cdot),v(\cdot)),s) f(\phi(s, t_0,x_0,\gamma(\cdot),v(\cdot)),\gamma(\cdot),v(\cdot)) \Big) ds \\
& < -\theta (t-t_0).
\end{align*}
Let $\tau_0 \in [t_0,t_0+ \delta_3]$ be such that
\begin{equation*}
h(\phi(\tau_0,t_0,x_0,\gamma(\cdot),v(\cdot))) = \max_{\tau \in [t_0,t_0+\delta_3]} h(\phi(\tau,t_0,x_0,\gamma(\cdot),v(\cdot))).
\end{equation*}
\textbf{Case 1.1:} If $\tau_0 \in (t_0,t_0+ \delta_3]$, then for $t=\tau_0$ we have
\begin{equation}
V(\phi(\tau_0,t_0,x_0,\gamma(\cdot),v(\cdot)),\tau_0) - V(x_0,t_0) < - \theta(\tau_0-t_0) < 0.
\end{equation}
Then by the dynamic programming argument of Lemma $1$ we have:
\begin{align*} V(x_0,t_0) \leq \sup_{v(\cdot)\in  \mathcal{V}_{[t_0,t_0+\delta_3]}} \big[ \max \big \{ &\max_{\tau \in [t_0,t_0+\delta_3]} h(\phi(\tau,t_0,x_0,\gamma(\cdot),v(\cdot))), \\&V(\phi(\tau_0,t_0,x_0,\gamma(\cdot),v(\cdot)),\tau_0) \big \} \big].
\end{align*}
We can choose $\hat{v}(\cdot) \in  \mathcal{V}_{[t_0,t_0+\delta_3]}$ such that
\begin{equation*} V(x_0,t_0) \leq \max \big \{ \max_{\tau \in [t_0,t_0+\delta_3]} h(\phi(\tau,t_0,x_0,\gamma(\cdot),v(\cdot))), V(\phi(\tau_0,t_0,x_0,\gamma(\cdot),v(\cdot)),\tau_0) \big \} + \epsilon,
\end{equation*}
and set $\epsilon < \frac{\theta}{2}(\tau_0-t_0)$.
Since $h(x)-V(x,t) \leq 0$ for all $(x,t) \in \mathbb{R}^n \times (0,T)$ we have that
\begin{equation*}\max_{\tau \in [t_0,t_0+\delta_3]} h(\phi(\tau,t_0,x_0,\gamma(\cdot),\hat{v}(\cdot))) = h(\phi(\tau_0,t_0,x_0,\gamma(\cdot),\hat{v}(\cdot))) \leq V(\phi(\tau_0,t_0,x_0,\gamma(\cdot),\hat{v}(\cdot)),\tau_0). \end{equation*}
Hence
\begin{equation*} V(x_0,t_0) \leq  V(\phi(\tau_0,t_0,x_0,\gamma(\cdot),\hat{v}(\cdot)),\tau_0) + \frac{\theta}{2}(\tau_0-t_0).\end{equation*}
Since $(10)$ holds for all $v(\cdot) \in \textit{V}_{[t_0,T]}$, it will also hold for $\hat{v}(\cdot)$, and hence the last argument establishes a contradiction. \\
\newline
\textbf{Case 1.2:} If $\tau_0=t_0$ then for $t=t_0 + \delta_3$ we have that for all $v(\cdot) \in \textit{V}_{[t_0,T]}$
\begin{equation*}
V(\phi(t_0 + \delta_3,t_0,x_0,\gamma(\cdot),v(\cdot)),t_0 + \delta_3) - V(x_0,t_0) < - \theta \delta_3 < 0.
\end{equation*}
Since by Lemma 1
\begin{align*} V(x_0,t_0) \leq \sup_{v(\cdot)\in  \mathcal{V}_{[t_0,t_0+\delta_3]}} \max \big \{ &\max_{\tau \in [t_0,t_0+\delta_3]} h(\phi(\tau,t_0,x_0,\gamma(\cdot),v(\cdot))), \\&V(\phi(t_0+ \delta_3,t_0,x_0,\gamma(\cdot),v(\cdot)),t_0+ \delta_3) \big \}, \end{align*}
then if
\begin{equation*} V(x_0,t_0) \leq \sup_{v(\cdot)\in  \mathcal{V}_{[t_0,t_0+\delta_3]}} V(\phi(t_0+ \delta_3,t_0,x_0,\gamma(\cdot),v(\cdot)),t_0+ \delta_3) , \end{equation*}
we can choose $\hat{v}(\cdot) \in  \mathcal{V}_{[t_0,t_0+\delta_3]}$ such that
\begin{equation*} V(x_0,t_0) \leq V(\phi(t_0+ \delta_3,t_0,x_0,\gamma(\cdot),\hat{v}(\cdot)),t_0+ \delta_3) + \frac{\theta \delta_3}{2}, \end{equation*}
which establishes a contradiction.\\
If \begin{equation*} V(x_0,t_0) \leq \sup_{v(\cdot)\in  \mathcal{V}_{[t_0,t_0+\delta_3]}}  \max_{\tau \in [t_0,t_0+\delta_3]} h(\phi(\tau,t_0,x_0,\gamma(\cdot),v(\cdot))), \end{equation*}
then we can choose $\hat{v}(\cdot) \in  \mathcal{V}_{[t_0,t_0+\delta_3]}$ such that
\begin{equation*} V(x_0,t_0) \leq  \max_{\tau \in [t_0,t_0+\delta_3]} h(\phi(\tau,t_0,x_0,\gamma(\cdot),\hat{v}(\cdot))) + \epsilon, \end{equation*}
or equivalently $V(x_0,t_0) \leq h(x_0) + \epsilon$, since $\tau_0=t_0$.
Based on our initial hypothesis that $h(x_0) < V(x_0,t_0)$, there exists a $\delta>0$ such that $h(x_0) - V(x_0,t_0) < -2\delta$. If we take $\epsilon<\delta$ we establish a contradiction.\\
\newline
\textbf{Part 2.}
Consider an arbitrary $(x_0,t_0) \in \mathbb{R}^n \times (0,T)$ and a smooth $W : \mathbb{R}^n \times (0,T) \rightarrow \mathbb{R}$ such that $V-W$ has a local minimum at $(x_0,t_0)$. Then, there exists $\delta_1 > 0$ such that for all $(x,t) \in \mathbb{R}^n \times (0,T)$ with $|x-x_0|^2 + (t-t_0)^2 < \delta_1$
\begin{equation*}
(V-W)(x_0,t_0) \leq (V-W)(x,t).
\end{equation*}
We would like to show that
\begin{equation*}
   \max \{ h(x_0)-V(x_0,t_0), \frac{\partial W}{\partial t}(x_0,t_0)+\sup_{v \in \textit{V}} \inf_{u \in \textit{U}} \frac{\partial W}{\partial x}(x_0,t_0) f(x_0,u,v) \} \leq 0.
\end{equation*}
Since $V(x,t) \geq h(x)$ it suffices to show that
$\frac{\partial W}{\partial t}(x_0,t_0)+\sup_{v \in \textit{V}} \inf_{u \in \textit{U}} \frac{\partial W}{\partial x}(x_0,t_0) f(x_0,u,v) \leq 0$.
This implies that for all $v \in \textit{V}$ there exists a $u \in \textit{U}$ such that
\begin{equation*}
\frac{\partial W}{\partial t}(x_0,t_0)+\frac{\partial W}{\partial x}(x_0,t_0)f(x_0,u,v) \leq 0.
\end{equation*}
For the sake of contradiction assume that there exists $\hat{v} \in \textit{V}$ such that for all $u \in \textit{U}$ there exists $\theta > 0$ such that
\begin{equation*}
\frac{\partial W}{\partial t}(x_0,t_0)+\frac{\partial W}{\partial x}(x_0,t_0)f(x_0,u,\hat{v}) > 2\theta > 0.
\end{equation*}
Since $W$ is smooth, there exists $\delta_2 \in (0,\delta_1)$  such that for all $(x,t) \in \mathbb{R}^n \times (0,T)$ with $|x-x_0|^2 + (t-t_0)^2 < \delta_2$
\begin{equation*}
\frac{\partial W}{\partial t}(x,t)+\frac{\partial W}{\partial x}(x,t)f(x,u,\hat{v}) > \theta > 0.
\end{equation*}
Hence, following \cite{Evans_Souganidis}, for $v(\cdot) \equiv \hat{v}$ and any $\gamma(\cdot) \in \Gamma_{[t_0,T]}$
\begin{equation*}
\frac{\partial W}{\partial t}(x,t)+\frac{\partial W}{\partial x}(x,t)f(x,\gamma(\cdot),v(\cdot)) > \theta > 0.
\end{equation*}
By continuity, there exists $\delta_3 > 0$  such that $|\phi(t,t_0,x_0,\gamma(\cdot),v(\cdot))-x_0|^2 + (t-t_0)^2 < \delta_2$ for all $t \in [t_0,t_0+\delta_3]$.
Therefore, for all $\gamma(\cdot) \in \Gamma_{[t_0,T]}$
\begin{align*}
V(\phi&(t_0 + \delta_3, t_0,x_0,\gamma(\cdot),v(\cdot)),t_0 + \delta_3)-V(x_0,t_0) \\&\geq W(\phi(t_0 + \delta_3, t_0,x_0,\gamma(\cdot),v(\cdot)),t_0 + \delta_3)-W(x_0,t_0) \\
& = \int_{t_0}^{t_0+\delta_3} \Big( \frac{\partial W}{\partial t}(\phi(t, t_0,x_0,\gamma(\cdot),v(\cdot)),t)\\&+ \frac{\partial W}{\partial x}(\phi(t, t_0,x_0,\gamma(\cdot),v(\cdot)),t) f(\phi(t, t_0,x_0,\gamma(\cdot),v(\cdot)),\gamma(\cdot),v(\cdot)) \Big) dt \\
& > \theta \delta_3.
\end{align*}
But by the dynamic programming argument of Lemma $1$ we can choose a $\hat{\gamma}(\cdot)\in  \Gamma_{[t_0,T]}$ such that
\begin{align*}
V(x_0,&t_0)  \geq \sup_{v(\cdot)\in  \mathcal{V}_{[t_0,t_0+\delta_3]}} \big[ \max \big \{ \max_{\tau \in [t_0,t_0+\delta_3]} h(\phi(\tau,t_0,x_0,\hat{\gamma}(\cdot),v(\cdot))), \\&V(\phi(t_0+ \delta_3,t_0,x_0,\hat{\gamma}(\cdot),v(\cdot)),t_0+ \delta_3) \big \} \big]  - \frac{\delta_3 \theta}{2} \\
& \geq  \max \big \{ \max_{\tau \in [t_0,t_0+\delta_3]} h(\phi(\tau,t_0,x_0,\hat{\gamma}(\cdot),v(\cdot))), V(\phi(t_0+ \delta_3,t_0,x_0,\hat{\gamma}(\cdot),v(\cdot)),t_0+ \delta_3) \big \} - \frac{\delta_3 \theta}{2} \\
& \geq V(\phi(t_0+ \delta_3,t_0,x_0,\hat{\gamma}(\cdot),v(\cdot)),t_0+ \delta_3) - \frac{\delta_3 \theta}{2}.
\end{align*}
The last statement establishes a contradiction, and completes the proof.
\end{proof}

\subsection{Variational inequality for $\widetilde{V}$}
Consider the value function $\widetilde{V}$ defined in the previous section. The following theorem proposes that $\widetilde{V}$ is the unique viscosity solution of another variational inequality.\\
\newline
\textbf{Theorem 2.} \textit{ $\widetilde{V}: \mathbb{R}^n \times [0,T] \rightarrow \mathbb{R}$ is the unique viscosity solution of the variational inequality }
\begin{equation}
    \max \Big \{ h(x)-\widetilde{V}(x,t), \frac{\partial \widetilde{V}}{\partial t}(x,t)+ \min \{0,\sup_{v \in \textit{V}} \inf_{u \in \textit{U}} \frac{\partial \widetilde{V}}{\partial x}(x,t) f(x,u,v)\} \Big \}=0,
\end{equation}
\textit{ with terminal condition $\widetilde{V}(x,T)= \max \{l(x),h(x) \}$.}
\begin{proof}
By Theorem 1, $\widetilde{V}(x,t)$ is the unique viscosity solution of $(4)$, subject to $\widetilde{V}(x,T)=\max\{l(x),h(x)\}$. If we let $\widetilde{H}(x,p) = \sup_{v \in \textit{V}} \inf_{\tilde{u} \in \widetilde{U}} p^T \tilde{f}(x,u,v)$ then, following the proof of Theorem 2 of \cite{Mitchell_Bayen2}, we have that
\begin{align*}
\widetilde{H}(x,p) &= \sup_{v \in \textit{V}} \inf_{\tilde{u} \in \widetilde{U}} p^T \tilde{f}(x,\tilde{u},v) \\
&= \sup_{v \in \textit{V}} \inf_{u \in \textit{U}} \inf_{\bar{u} \in \bar{U}} p^T (\bar{u}f(x,u,v)) \\
&= \inf_{\bar{u} \in \bar{U}} \bar{u} \sup_{v \in \textit{V}} \inf_{u \in \textit{U}}  p^T f(x,u,v) \\
&= \min_{\bar{u} \in \bar{U}} \bar{u} \sup_{v \in \textit{V}} \inf_{u \in \textit{U}}  p^T f(x,u,v) \\
&= \min\{0,H(x,p)\}.
\end{align*}
Consequently, the two variational inequalities $(4)$ and $(11)$ are equivalent, and so $\widetilde{V}(x,t)$ is the viscosity solution of $(11)$.
\end{proof}
Since the solution to $(11)$ is unique \cite{Georgiou}, one could easily show that
\begin{equation*}
  \widetilde{V}(x,t)=\inf_{\gamma(\cdot) \in  \Gamma_{[t,T]}}  \sup_{v(\cdot) \in  \mathcal{V}_{[t,T]}} \min_{\tau_1 \in [t,T]} \max \{l(\phi(\tau_1,t,x,u(\cdot),v(\cdot))) , \max_{\tau_2 \in [t,\tau_1]} h(\phi(\tau_2,t,x,u(\cdot),v(\cdot))) \}.
\end{equation*}

\section{Case study: Collision Avoidance in Air Traffic Management}
To illustrate the approach described in the previous sections, we consider a problem from the air traffic management area. The increase in air traffic is bound to lead to further en-route delays and potentially safety problems in the immediate future \cite{SESAR}, \cite{NextGen}. A major difficulty with accommodating this expected increase in air traffic is uncertainty about the future evolution of flights. Therefore, the CATS research project has proposed a novel concept of operations, which aims to increase punctuality and safety during the flight. This concept is mainly based on imposing spatial and temporal constraints at different parts of the flight plan of each aircraft. These 4D constraints are known as Target Windows (TW) \cite{CATS1}, and represent the commitment from each actor (air traffic controllers, airports, airlines, air navigation service providers) to deliver a particular aircraft within the TW constraint. This commitment is known as the Contract of Objectives (CoO) \cite{CATS}, and can be viewed as a first step towards the implementation of the Reference Business Trajectory envisioned by the SESAR joint undertaking \cite{SESAR}.

In this section we follow the approach proposed in \cite{CDCpaper} to code the TW constraints, and use the reach-avoid formulation of Sections II and III, to investigate collision avoidance in the presence of TW constraints. For this purpose we consider a two-aircraft scenario, where each aircraft should respect its TW constraints, while avoiding conflict with other aircraft.

\subsection{Aircraft model}
Each aircraft $j=1,...,N$ is assumed to have a predetermined flight plan, which comprises a series of way points $ O_{(i,j)} = \left[x_{(i,j)}~y_{(i,j)}~z_{(i,j)}\right]^T \in \mathbb{R}_+^3$, where $i=1,...,M_{j}$. The angle $\Psi_{(i,j)}$ that each segment forms with the $x$ axis and the flight path angle $\Gamma_{(i,j)}$ that it forms with the horizontal plane are shown in Fig. 1. The discrete state $i$ stores the segment of the flight plan that the aircraft is currently in, and for $i=1,...,M_{j}-1$ we can define
\begin{equation*}
\begin{array}{rl}
%d_i &= \vectornorm{ \begin{bmatrix} x_{i+1} \\ y_{i+1} \end{bmatrix} - \begin{bmatrix} x_i \\ y_i \end{bmatrix} }, \\
\Psi_{(i,j)} &= \tan^{-1}\left( \frac{y_{(i+1,j)}-y_{(i,j)}}{x_{(i+1,j)}-x_{(i,j)}} \right), \\
\Gamma_{(i,j)} &= \tan^{-1}\left( \frac{z_{(i+1,j)}-z_{(i,j)}}{d_{(i,j)}} \right),
\end{array}
\end{equation*}
where $d_{(i,j)} = \sqrt{(x_{(i+1,j)}-x_{(i,j)})^2 + (y_{(i+1,j)}-y_{(i,j)})^2}$ is the length of the projection of its segment on the horizontal plane. Assume perfect lateral tracking and set $s \in \mathbb{R}_+$ denote the the part of each segment covered on the horizontal plane (see Fig. 1).
Based on our assumption that each aircraft has constant heading angle $\Psi_{(i,j)}$ at each segment, its $x$ and $y$ coordinates can be computed by:
\begin{equation*}
\begin{bmatrix} x_{(i,j)}(s_j) \\ y_{(i,j)}(s_j) \end{bmatrix} = \begin{bmatrix} x_{(i,j)} \\ y_{(i,j)} \end{bmatrix} + \begin{bmatrix} \cos{\Psi_{(i,j)}} \\ \sin{\Psi_{(i,j)}} \end{bmatrix} s_j .
\end{equation*}
%\begin{figure}[!htb]
%\centering
%\begin{minipage}[b]{0.45\linewidth}
%\centering
%\includegraphics[scale=0.45]{f1}
%\caption{Flight plan projection on the horizontal plane}\label{fig:f1}
%\end{minipage}
%\hspace{0.45cm}
%\begin{minipage}[b]{0.45\linewidth}
%\centering
%\includegraphics[scale=0.45]{f2}
%\caption{Flight plan projection on the $z$-$s$ plane}\label{fig:f2}
%\end{minipage}
%\end{figure}

\begin{figure}[htp]
\centering
\subfigure[Flight plan projection on the horizontal plane]{
\includegraphics[scale=0.45]{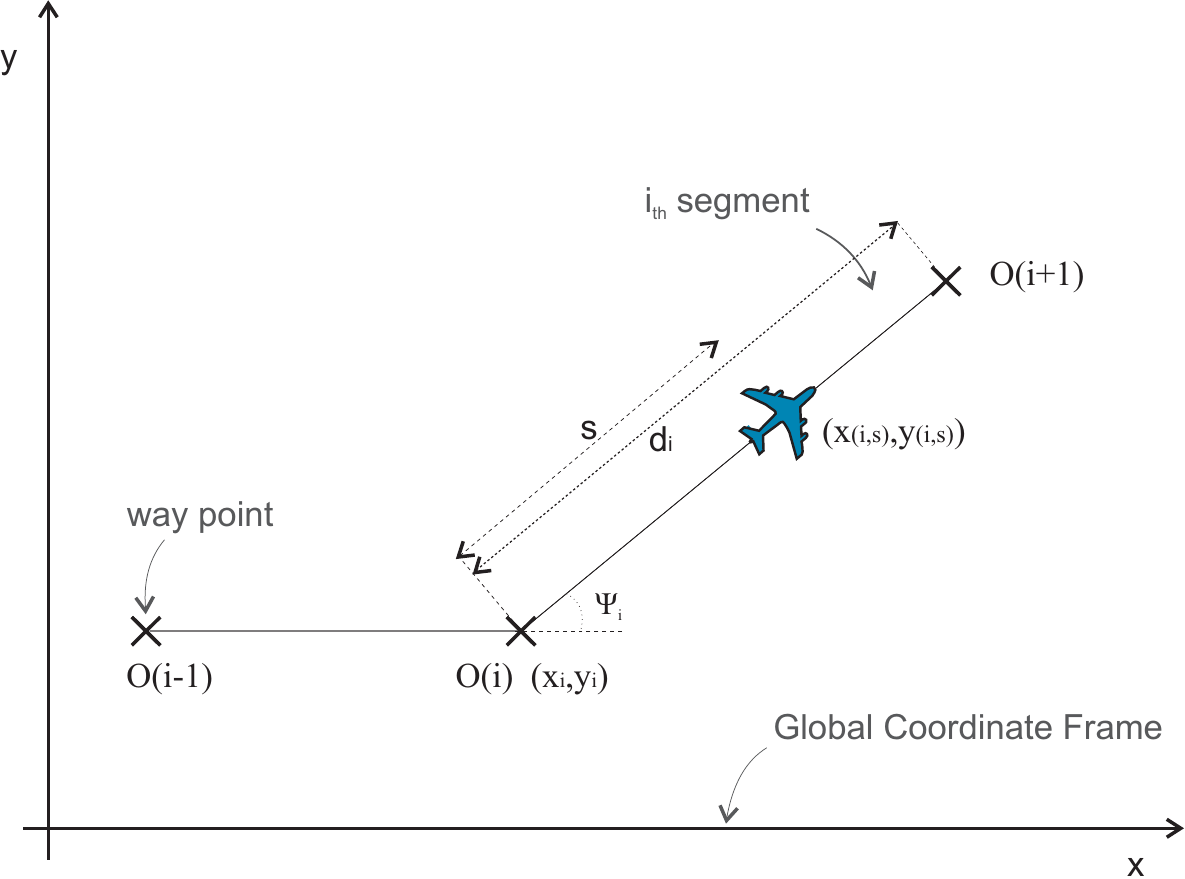}
\label{fig:f1}
}
\subfigure[Flight plan projection on the $z$-$s$ plane]{
\includegraphics[scale=0.45]{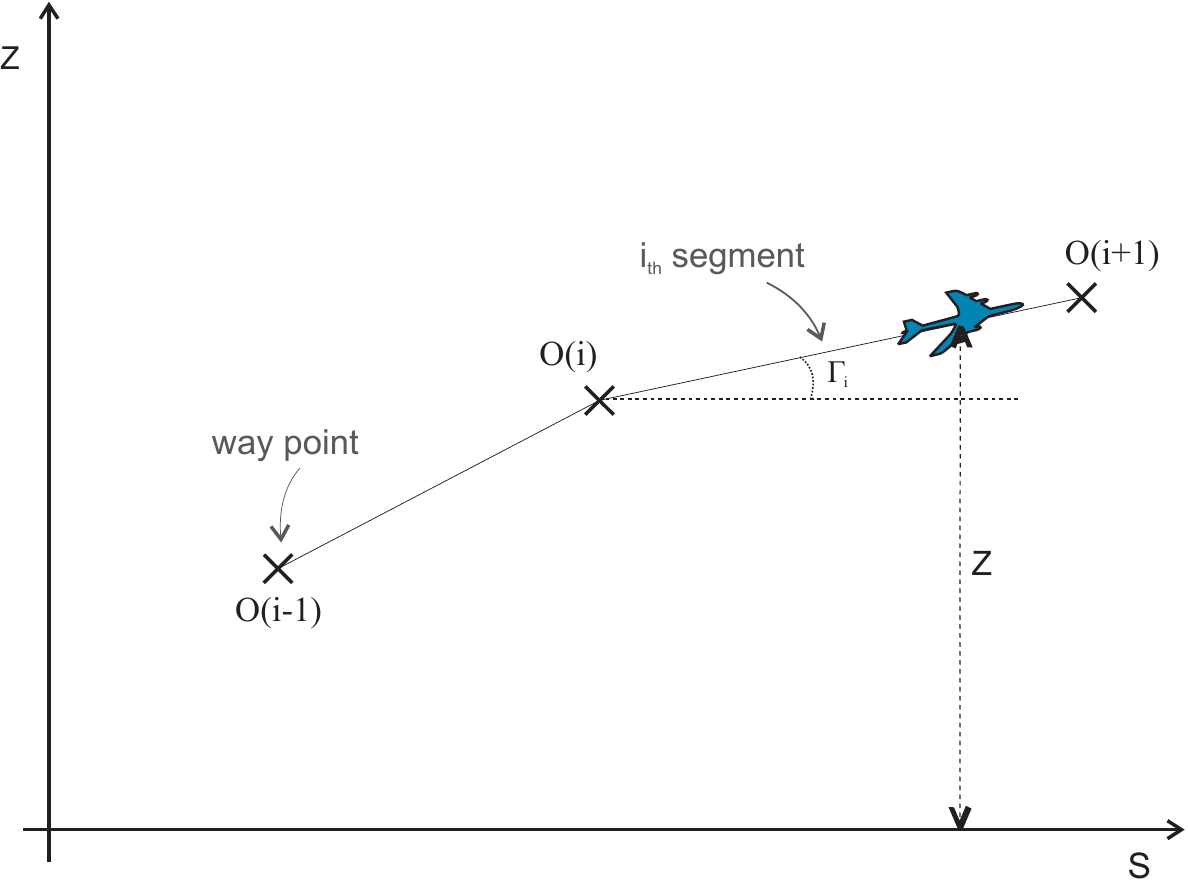}
\label{fig:f2}
}
\label{fig:subfigureExample}
\caption[Flight plan projections.]{Flight plan projections}
\end{figure}
To approximate accurately the physical model, the flight path angle $\gamma^p_{j}$, which is the angle that the aircraft forms with the horizontal plane, is a control input fixed according to the angle $\Gamma_{(i,j)}$ that the segment forms with the horizontal plane. If $\Gamma_{(i,j)}=0$ the aircraft will be cruising $(\gamma^p_{j}=0)$ at that segment, whereas if it is positive or negative it will be climbing $(\gamma^p_{j} \in [0,\overline{\gamma}^p_j])$ or descending $(\gamma^p_{j} \in [-\overline{\gamma}^p_j,0])$ respectively.

The speed of each aircraft apart from its type depends also on the altitude $z$. At each flight level there is a nominal airspeed that aircraft tend to track, giving rise to a function $g(z_{j},\gamma^p_{j})$. The dependence on the flight path angle indicates the discrete mode i.e. cruise, climb, descent, that an aircraft could be. For our simulations, we have assumed that at every level the airspeed could vary within $10\%$ of the nominal one; this is restricted by the control input $b_{j} \in [-1,1]$. Figure 2 shows the speed-altitude profiles of the A320 (the simulated aircraft), based on the BADA database \cite{bada}, for the different phases of flight. These curves have been computed by linear interpolation between the predetermined '$.$' points.

%\begin{figure}[htp]
%%\centering
%\begin{minipage}[b]{0.45\linewidth}
%\includegraphics[scale=0.45]{climb}
%\caption{CLIMB Speed-Altitude profile}\label{fig:climb}
%\end{minipage}
%\hspace{0.45cm}
%%\centering
%\begin{minipage}[b]{0.45\linewidth}
%\includegraphics[scale=0.45]{descent}
%\caption{DESCENT Speed-Altitude profile}\label{fig:descent}
%\end{minipage}
%\end{figure}
%\begin{figure}[h]
%\centering
%\includegraphics[scale=0.45]{cruise}
%\caption{CRUISE Speed-Altitude profile}\label{fig:cruise}
%\end{figure}

\begin{figure}[htp]
\centering
\subfigure[CLIMB Speed-Altitude profile]{
\includegraphics[scale=0.45]{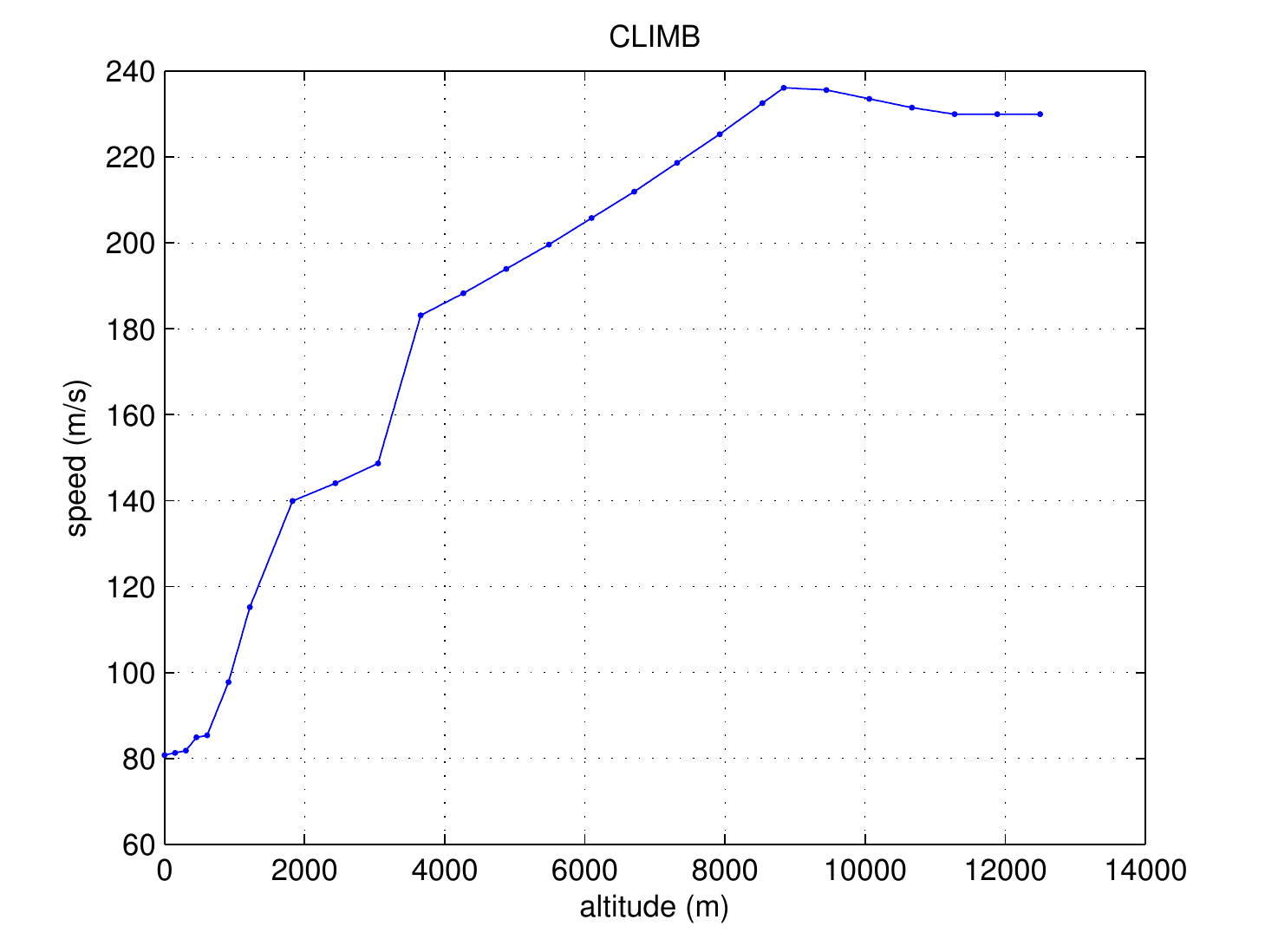}
\label{fig:climb}
}
\subfigure[DESCENT Speed-Altitude profile]{
\includegraphics[scale=0.45]{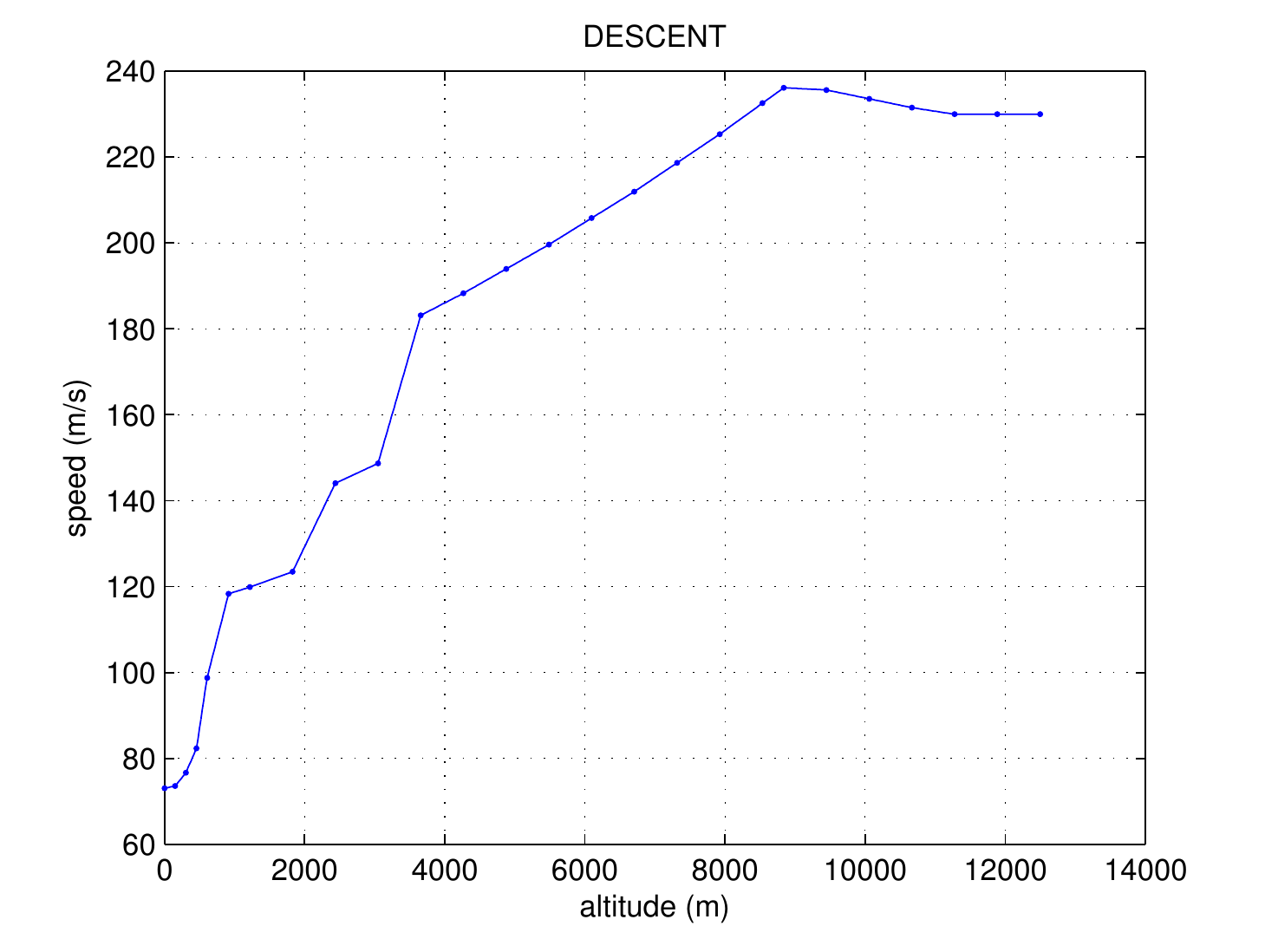}
\label{fig:descent}
}
\subfigure[CRUISE Speed-Altitude profile]{
\includegraphics[scale=0.45]{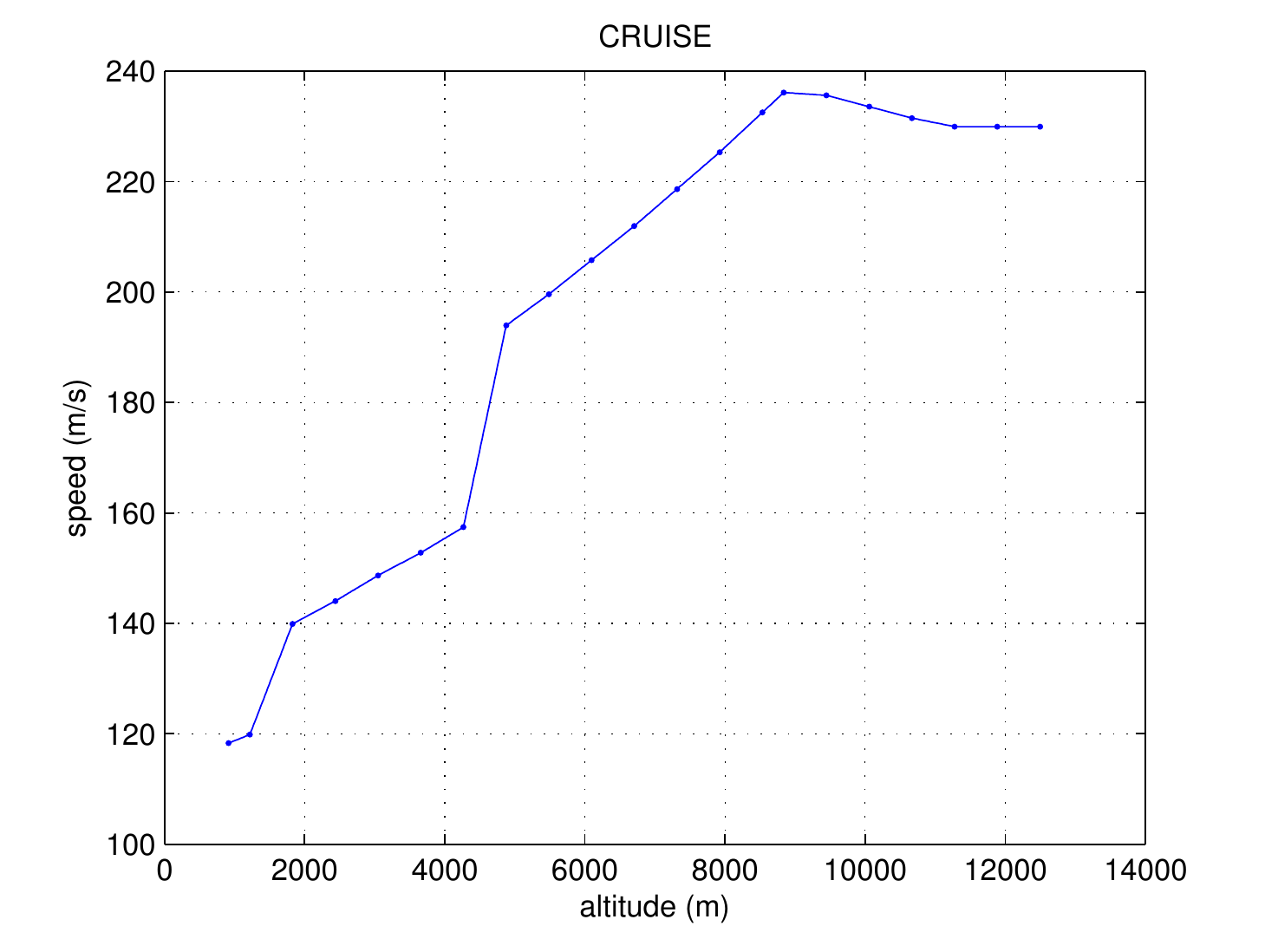}
\label{fig:cruise}
}
\label{fig:subfigureExample}
\caption[Speed-Altitude profiles for the different phases of flight.]{Speed-Altitude profiles for the different phases of flight}
\end{figure}

Most of the reachability numerical methods are based on gridding the state space, so the memory and time necessary for the computation grow exponentially in the state dimension. Therefore, using a full, five- or six-state, point mass model of the aircraft, like the one described in \cite{Yannis}, would be computationally expensive to analyze using the existing computational tools.
In \cite{CDCpaper}, the full point mass model for the aircraft, was abstracted to a simplified one to make the reachability computation tractable. Motivated by the fact that aircraft track laterally very well, it was assumed that the heading angle remains constant at each segment. The dynamics of each aircraft are modeled by a hybrid automaton $H_j =(X_j,Q_j,Init_j,f_j,Dom_j,G_j,R_j)$ (in the notation of \cite{Gao_Lygeros}), with:
\begin{itemize}
\item continuous states $x_j=\left[s_j~z_j~t\right]^T \in \mathbb{R}_+^3 = X_j$.
\item discrete states $i \in \{0,...,M_j-1\} = Q_j$.
\item initial states $Init_j = \{(i,s_j,z_j,t) ~|~ i=0,s_j=0,z_j=z_{j_0}\}$.
\item control inputs $u_j = \left[b_j~\gamma^p_j\right]^T\in \left[-1,1\right] \times \left[-\overline{\gamma_j}^p,\overline{\gamma_j}^p\right]$.
\item disturbance inputs $v = \left[w_x~w_y~w_z\right]^T\in\mathbb{R}^3$
\item vector field $f_j:Q_j \times X_j \times U_j \times V \to X_j $.
\end{itemize}
\[ f_j(i,s_j,z_j,t)=\begin{bmatrix} \dot{s}_j \\ \dot{z}_j \\ \dot{t} \end{bmatrix} =
\begin{bmatrix} (1+0.1b_j)g(z_j,\gamma^p_j)\cos{\gamma^p_j} + w_x\cos{\Psi_{(i,j)}}+w_y\sin{\Psi_{(i,j)}} \\
(1+0.1b_j)g(z_j,\gamma^p_j)\sin{\gamma^p_j} + w_z \\ 1 \end{bmatrix} \]
\begin{itemize}
\item domain $Dom_j=\{(i,s_j,z_j,t)~|~s_j\leq d_{(i,j)}\}$.
\item guards $G_j(i,i+1)=\{(s_j,z_j,t)~|~s_j > d_{(i,j)}\}$.
\item reset map $R_j(i,i+1,s_j,z_j,t)=\{(0,z_j,t)\}$.
\end{itemize}

Apart from $s_j$, the other two continuous states are the altitude $z_j$, and the time $t$.
The last equation was included in order to track the TW temporal constraints. As stated above, $\gamma^p_j$ is the flight path angle and $w$ is the wind speed, which acts as a bounded disturbance
with $-\overline{w}\leq w\leq \overline{w}$, and for our simulations we used $\overline{w} = 12m/s$. Since the flight path angle $\gamma^p_{j}$ does not exceed 5$^\circ$, for simplification we can assume that $\sin \gamma^p_{j} \approx \gamma^p_{j}$ and $\cos \gamma^p_{j} \approx 1$.

\subsection{Reach-Avoid problem formulation}
Target Windows represent spatial and temporal constraints that aircraft should respect. Following \cite{CATS1}, we assume that TW are located on the surface area between two air traffic control sectors. Based on the structure of those sectors, the TW are either adjacent or superimposed (Fig. 3), and for simplicity we assume that there is a way point centered in the middle of each TW.

Our objective is to compute the set $I$ of all initial states at time $t$ for which there exists a non-anticipative control strategy $\gamma$, that despite the wind input $v$ can lead the aircraft $j$ inside the TW constraint set at least once within its time and space window, while avoiding conflict with the other aircraft.
In air traffic, conflict refers to the loss of minimum separation between two aircraft. Each aircraft is surrounded by a protected zone, which is generally thought of as a cylinder of radius 5nmi and height 2000ft centered at the aircraft. If this zone is violated by another aircraft, then a conflict is said to have occurred.
To achieve this goal, we adopt another simplification introduced in \cite{CDCpaper}; we eliminate time from the state equations, and perform a two-stage calculation.
\begin{figure}[htp]
\centering
\includegraphics[scale=0.95]{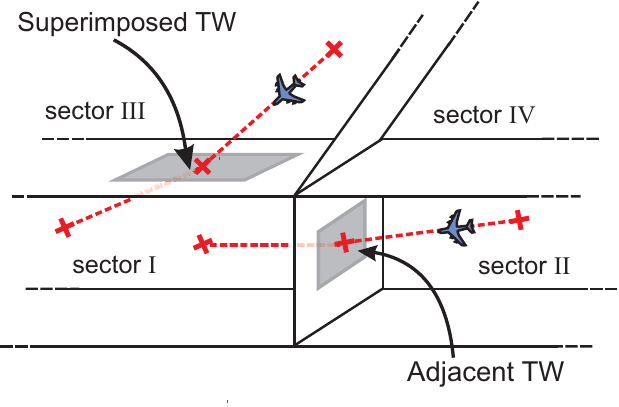}
\caption{Superimposed and Adjacent TWs}\label{fig:TWs_new}
\end{figure}

We define the spatial constraints of a TW centered at the way point $i$ as $\widetilde{R}_j = (d_{(i,j)},[z_{(i,j)}+\underline{z}_{(i,j)},z_{(i,j)}+\overline{z}_{(i,j)}])$ if the TW is adjacent, and $\widetilde{R}_j = ([d_{(i,j)}+\underline{s}_{(i,j)}, d_{(i,j)}+\overline{s}_{(i,j)}],z_{(i,j)})$ if the TW is superimposed. Let also $[\underline{t}_j, \overline{t}_j]$ denote the time window of $\widetilde{R}_j$. Then $I$ could be computed as:\\
\textbf{Stage 1:} Compute for each aircraft $j$ the set $R_j$ of states $x_j$ at time $\underline{t}_j$ (beginning of target window) from which there exist a control trajectory that despite the wind can lead the aircraft inside $\widetilde{R}_j$ at least once within the time interval $[\underline{t}_j, \overline{t}_j]$. But this set is the $\widetilde{RA}_j(\underline{t},\widetilde{R}_j,A_j)$ set, which was shown in Section III to be the zero sublevel set of $\widetilde{V}$, which is the solution to the following partial differential equation
\begin{equation*}
    \max \{ h_j(x_j)-\widetilde{V}(x_j,t), \frac{\partial \widetilde{V}}{\partial t}(x_j,t)+ \min \{0,\sup_{v \in \textit{V}} \inf_{u_j \in \textit{U}_j} \frac{\partial \widetilde{V}}{\partial x_j}(x_j,t) f_j(x_j,u_j,v)\} \}=0.
\end{equation*}
The terminal condition $\widetilde{V}(x_j,\overline{t}_j)=l(x_j)$ was chosen to be the signed distance to the set $\widetilde{R}_j^c$, and the avoid set $A_j$ is characterized by $h_j(x_j)$. This function represents the area where a conflict might occur, and it is computed online by performing conflict detection (see Appendix C). \\
\textbf{Stage 2:} Compute the set $I$ of all states that start at time $t \leq \underline{t}_j$ and for every wind can reach the set $R_j$ at time $\underline{t}_j$, while avoiding conflict with other aircraft. Based on the analysis of Section II, this is the $RA_j(t,R_j,A_j)$ set, that can be computed by solving
\begin{equation*}
    \max \{ h_j(x_j)-V(x_j,t), \frac{\partial V}{\partial t}(x_j,t)+\sup_{v \in \textit{V}} \inf_{u_j \in \textit{U}_j} \frac{\partial V}{\partial x_j}(x_j,t) f(x_j,u_j,v) \}=0,
\end{equation*}
with terminal condition $V(x_j,\underline{t}_j)=\max\{\widetilde{V}(x_j,\underline{t}_j),h_j(x_j)\}$. The set $R_j$ is defined as
\begin{equation*}
    R_j = \{x_j \in \mathbb{R}^n |~ \widetilde{V}(x_j,\underline{t}_j)) \leq 0\},
\end{equation*}
whereas $A_j$ depends once again on the obstacle function $h_j(x_j)$. 
\begin{figure}[htp]
\centering
\subfigure[Flight plans for the two aircraft case]{
\includegraphics[scale=0.45]{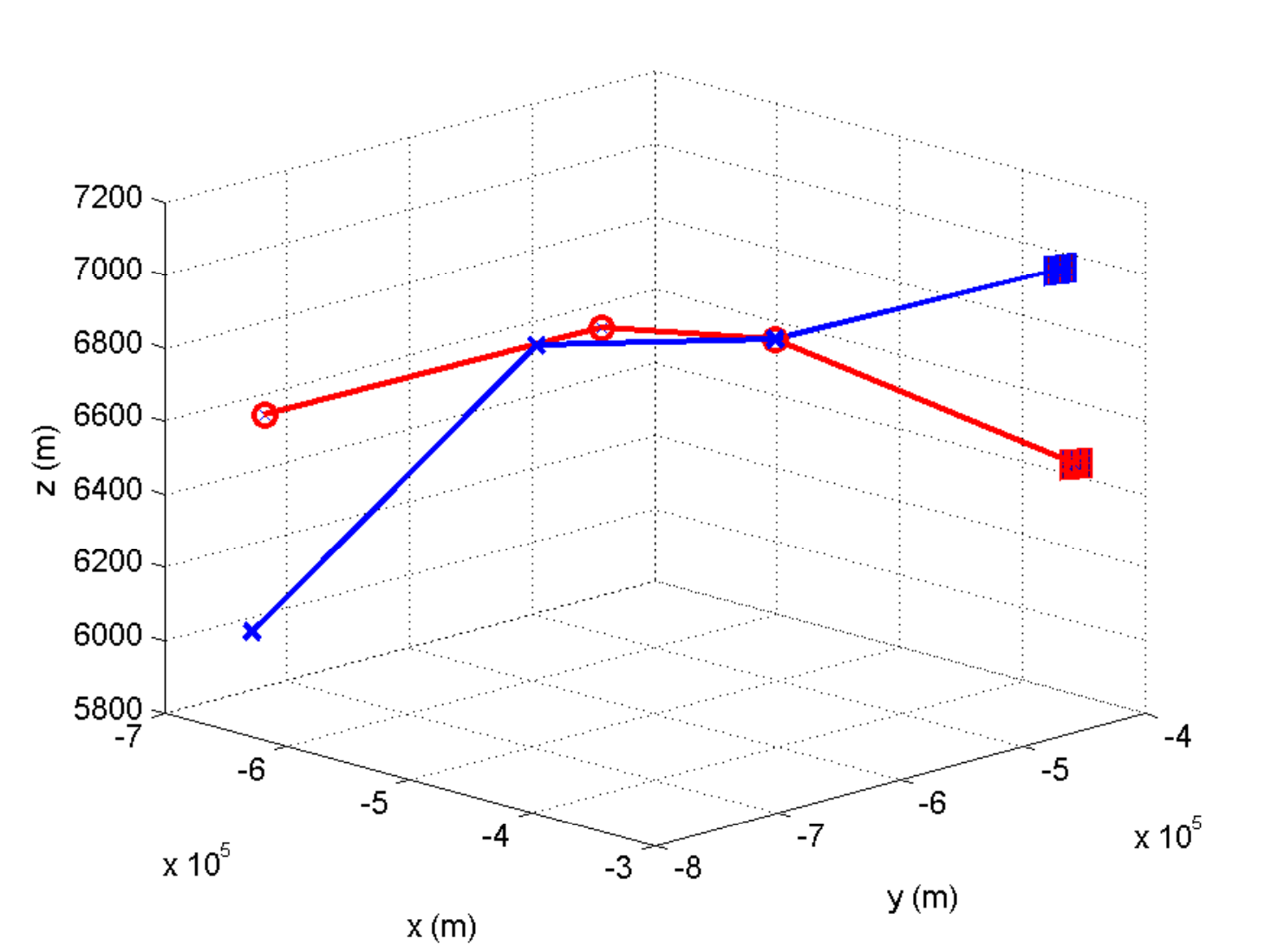}
\label{fig:flightPlansTW}
}
\subfigure[Flight plan projection for the two aircraft case]{
\includegraphics[scale=0.45]{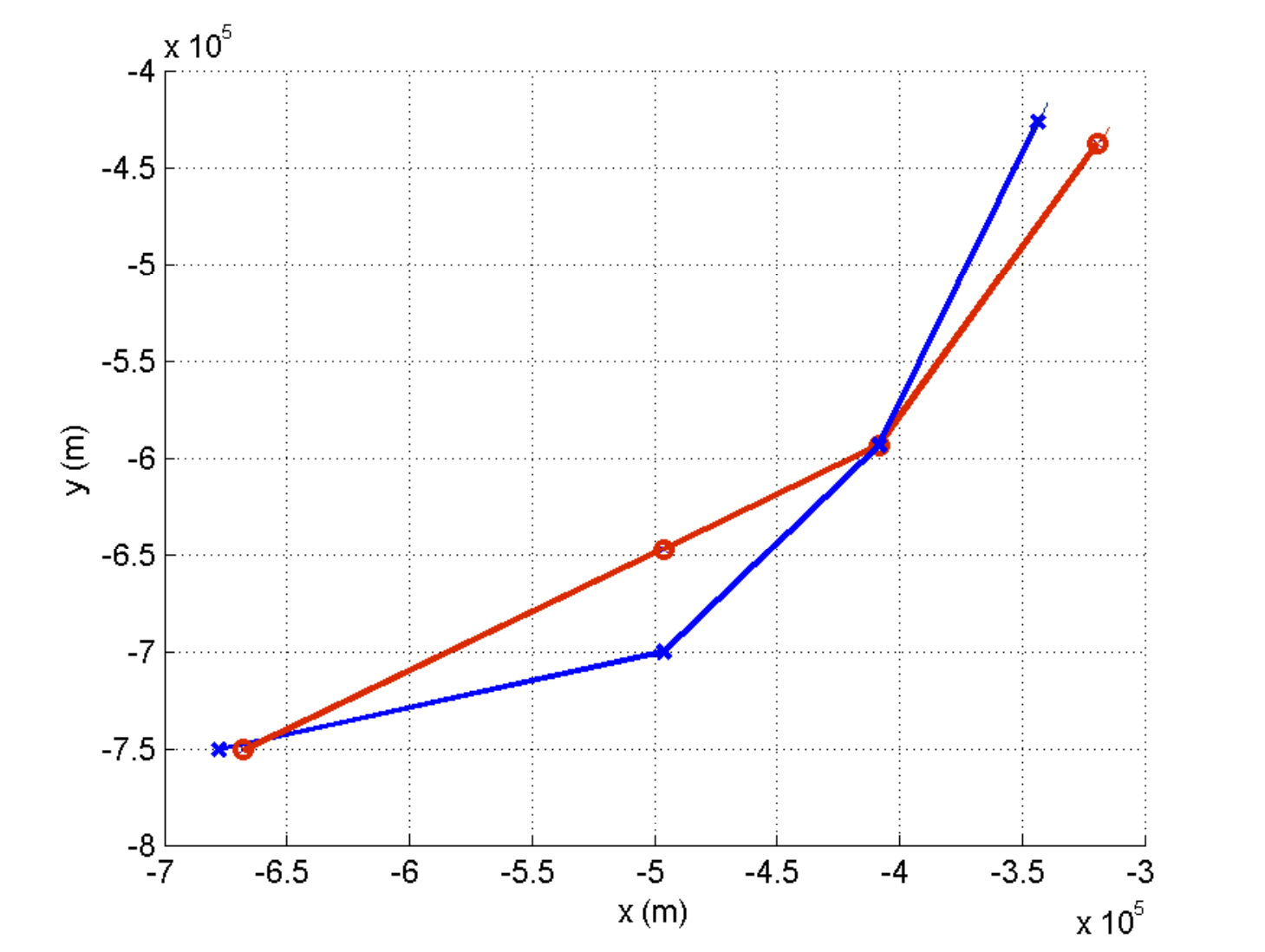}
\label{fig:fp2}
}
\label{fig:subfigureExample}
\caption[Flight plans for the two aircraft scenario]{Flight plans for the two aircraft scenario}
\end{figure}
The simulations for each aircraft are running in parallel, so at every instance $t$, we have full knowledge of the backward reachable sets of each aircraft. Based on that, Algorithm 1 of Appendix C describes the implemented steps for the Reach-Avoid computation.
\subsection{Simulation Results}
Consider now the case where we have two aircraft each one with a TW, whose flight plans intersect, and they enter the same air traffic sector with a $30sec$ difference. Fig. 4a depicts the two flight plans and Fig. 4b the projection of the flight plans on the horizontal plane. The Target Windows are centered at the last way point of each flight plan. The result of the two-stage backward reachability computation with TW as terminal sets is depicted in Fig. 5a. The tubes at this figure include all the states that each aircraft could be, and reach its TW. We should also note that the tubes are the union of the corresponding sets. These sets at a specific time instance, would include all the states that could start at that time and reach the TW at the end of the horizon. Fig. 5b is the projection of these tubes on the horizontal plane. As it was expected, the $x$-$y$ projection coincides with the projection of the flight plans on the horizontal plane. This is reasonable, since in the hybrid model we assumed constant heading angle at each segment. Moreover, based on the speed-altitude profiles, aircraft fly faster at higher levels, so at those altitudes there are more states that can reach the target.
%\begin{figure}[ht]
%\begin{minipage}[b]{0.45\linewidth}
%\centering
%\includegraphics[scale=0.45]{aircraft1J}
%\caption{Intersection of the two backward reachable tubes}\label{fig:aircraft1J}
%\end{minipage}
%\hspace{0.45cm}
%\begin{minipage}[b]{0.45\linewidth}
%\centering
%\includegraphics[scale=0.45]{aircraft2J}
%\caption{Projection of the two backward reachable tubes}\label{fig:aircraft2}
%\end{minipage}
%\end{figure}

\begin{figure}[htp]
\centering
\subfigure[Intersection of the two backward reachable tubes]{
\includegraphics[scale=0.45]{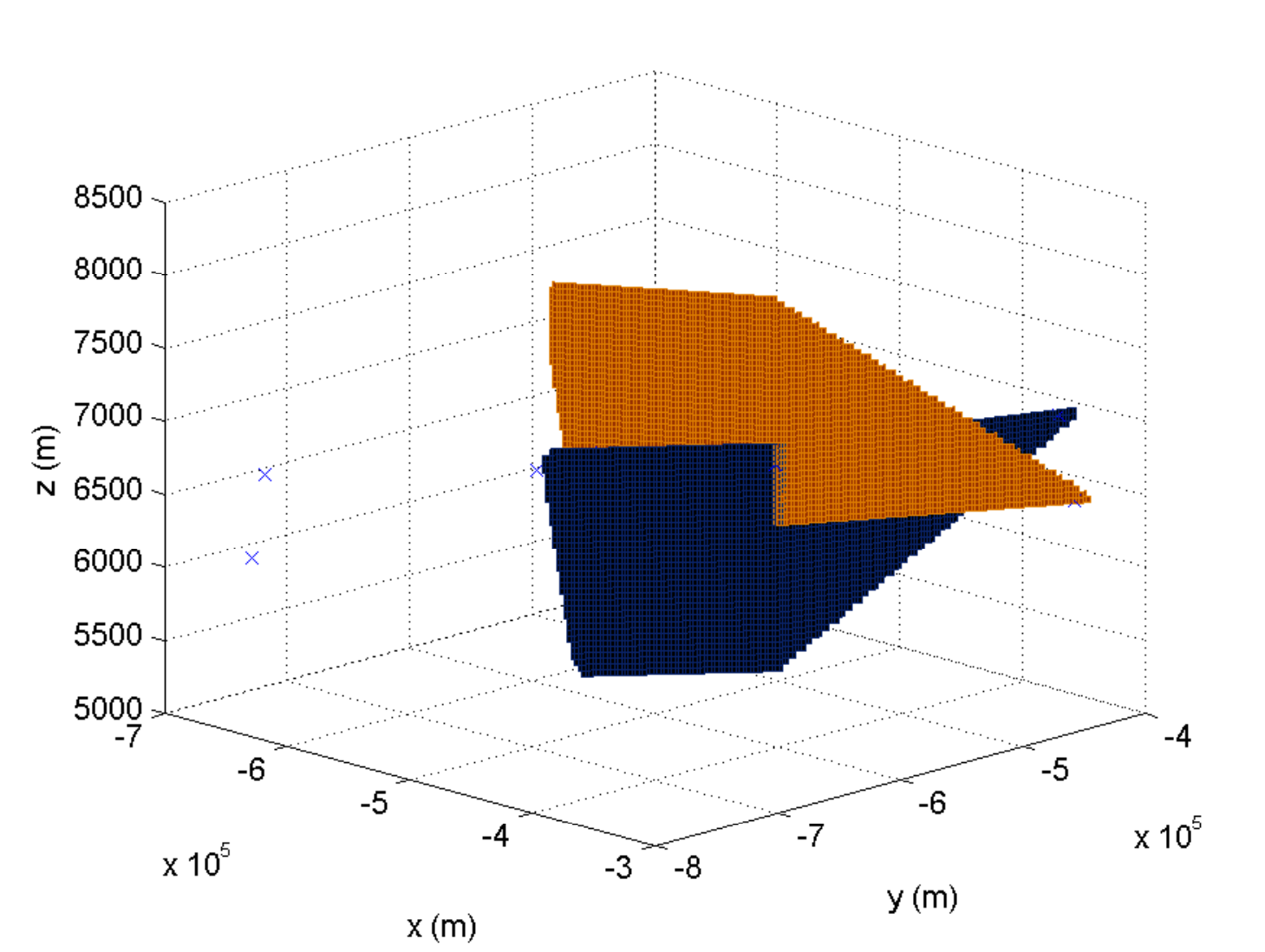}
\label{fig:aircraft1J}
}
\subfigure[Projection of the two backward reachable tubes]{
\includegraphics[scale=0.45]{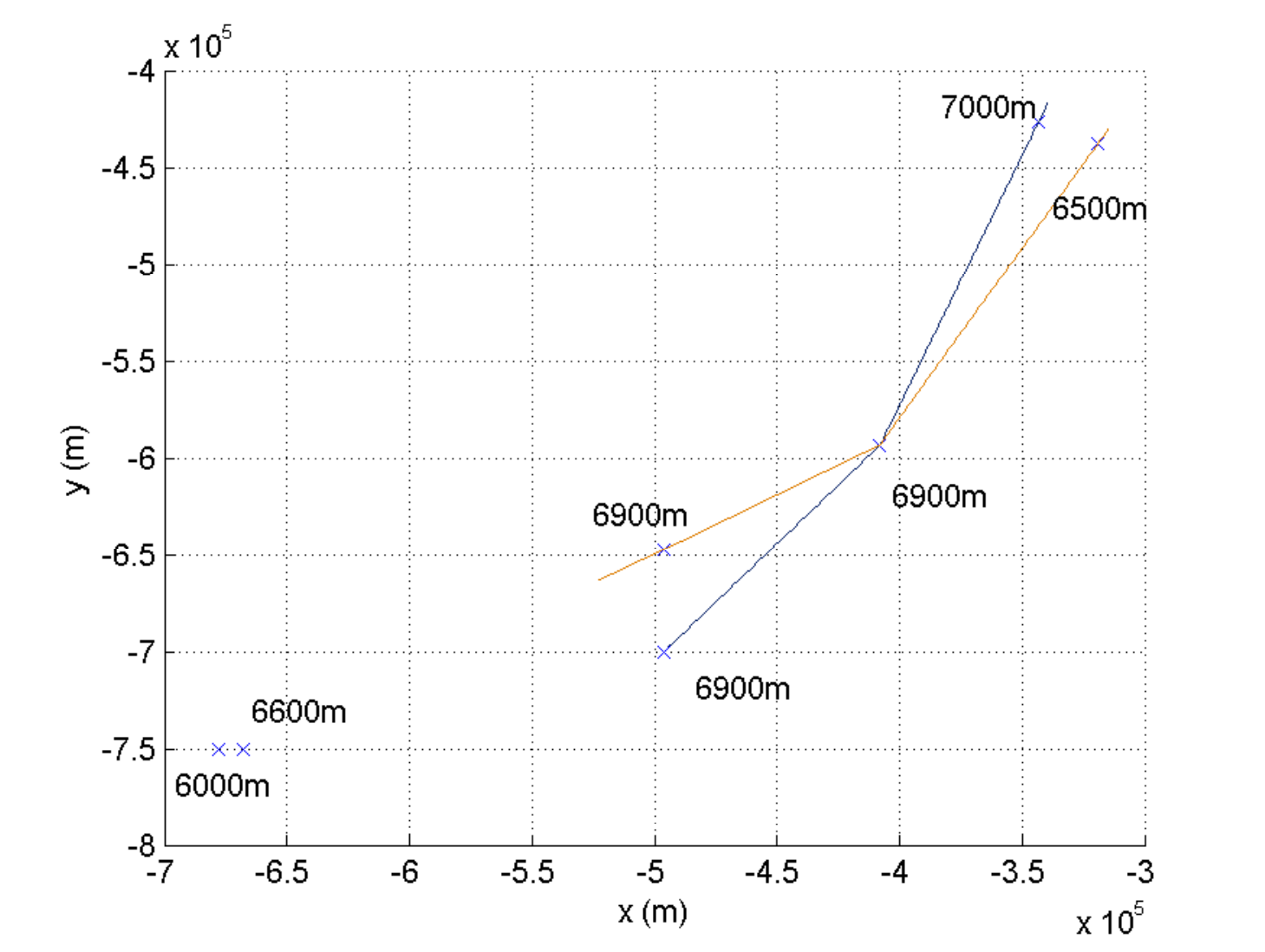}
\label{fig:aircraft2}
}
\label{fig:subfigureExample}
\caption[Backward reachable tubes for the two aircraft scenario]{Backward reachable tubes for the two aircraft scenario}
\end{figure}
We can repeat the previous computation, but now checking at every time if the sets, in the sense described before, satisfy the minimum separation standards. That way, the time and the points of each set where a conflict might occur, can be detected. The result of this calculation is illustrated in Fig. 6a. The "hole" that is now around the intersection area of Fig. 5a represents the area where the two aircraft might be in conflict.
%\begin{figure}[ht]
%\begin{minipage}[b]{0.45\linewidth}
%\centering
%\includegraphics[scale=0.45]{conf_area}
%\caption{Conflict detection zone}\label{fig:conf_area}
%\end{minipage}
%\hspace{0.45cm}
%\begin{minipage}[b]{0.45\linewidth}
%\centering
%\includegraphics[scale=0.45]{reachAvoidJ}
%\caption{Reach-Avoid tubes }\label{fig:reachAvoidJ}
%\end{minipage}
%\end{figure}

\begin{figure}[htp]
\centering
\subfigure[Conflict detection zone]{
\includegraphics[scale=0.45]{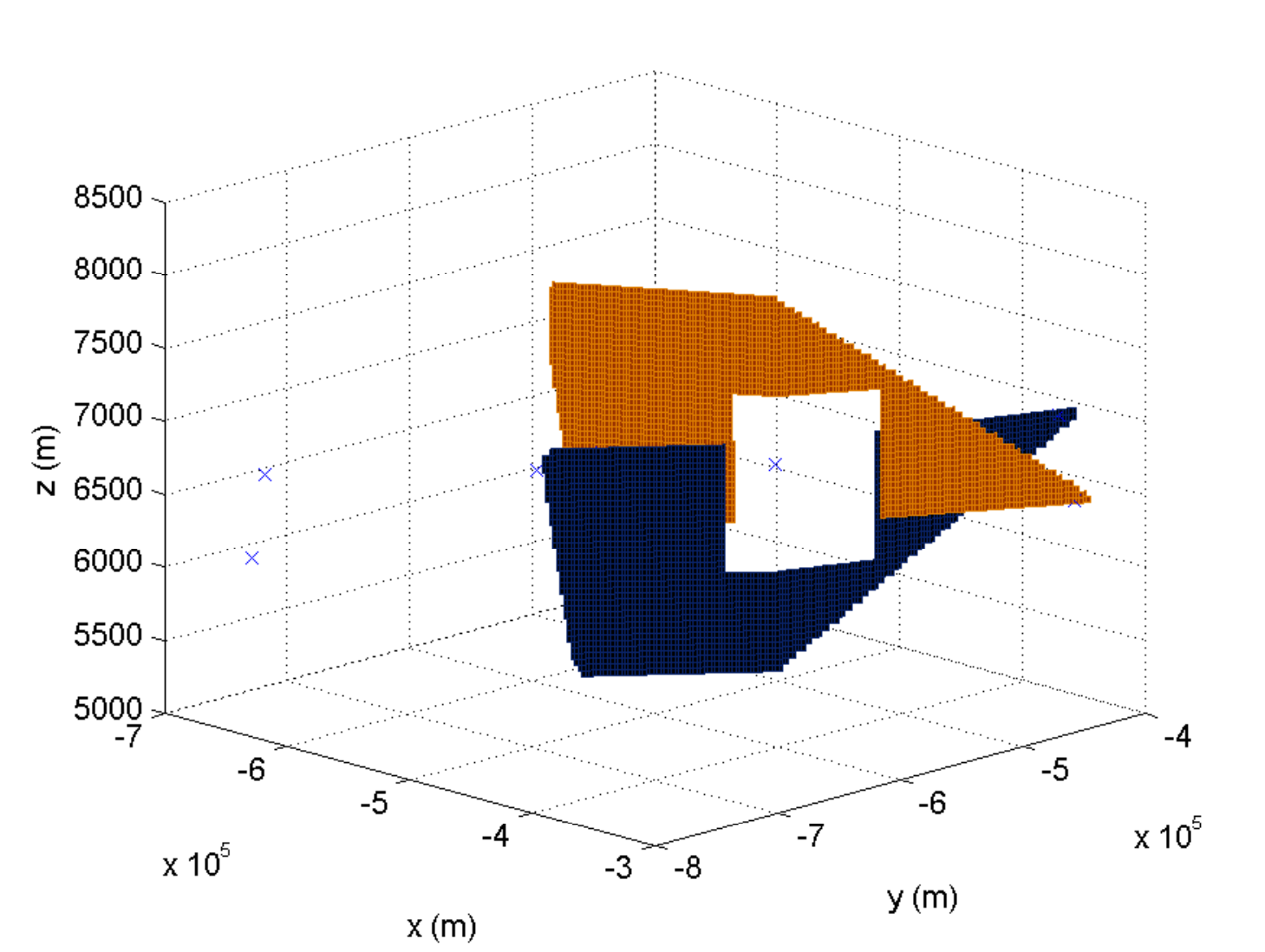}
\label{fig:conf_area}
}
\subfigure[Reach-Avoid tubes]{
\includegraphics[scale=0.45]{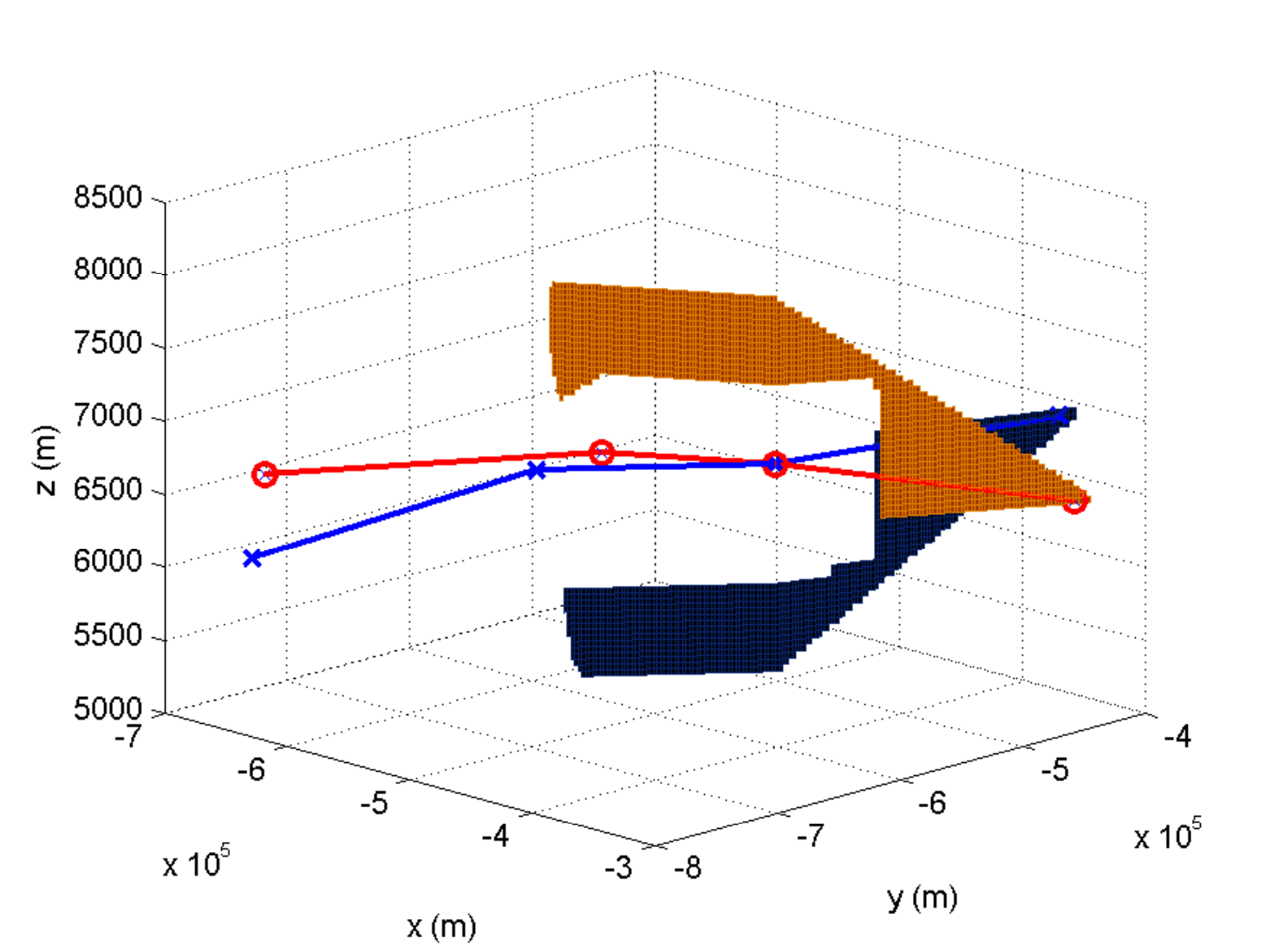}
\label{fig:reachAvoidJ}
}
\label{fig:subfigureExample}
\caption[Conflict detection and reach-avoid computation]{Conflict detection and reach-avoid computation}
\end{figure}
Now that we managed to perform conflict detection, we are in a position to compute at every instance the obstacle function $h_j(x_j)$. Since the conflict does not occur within the time interval of the TW, the set of the initial states that an aircraft $j$ could start and reach the set $R_j$ at time $\underline{t}_j$, while avoiding conflict with the other aircraft, should be computed. Once the aircraft hits $R_2$, it can also reach the TW within its time constraints. To obtain the solution to this reach-avoid problem, the variational inequality $(4)$ should be solved. If the conflict had occurred in $[\underline{t}_j,\overline{t}_j]$ equation $(11)$ should be solved instead. One could either use numerical methods developed by \cite{Georgiou}, or the Level Set Method Toolbox \cite{Mitchell_Tomlin}, whose authors propose a way to code obstacles on the value function. The latter was used in this paper, and the obstacle function $h_j(x_j)$ was dynamically determined, since at every time it is the result of the conflict detection. Fig. 6b shows the reach-avoid tubes at $t=25min$. As it was expected, the set of states that could reach the target while avoiding conflict with the other aircraft, does not include the conflict zone of Fig. 6a, but also some more states that would end up in this zone. This is a "centralized" solution, since for the safe sets shown in Fig. 6b, there exists a combination of inputs for the two aircraft that could satisfy all constraints, i.e. reaching their TW while avoiding conflict.
It should be also noted that all simulations where performed on the same grid, but were running in parallel. Hence he have two 2D computations (one for each aircraft) instead of a 4D one, and $h_j(x_j)$ was determined at each step by comparing the obtained sets as described in Algorithm 1.

\section{Concluding Remarks}
A new framework of solving nonlinear systems with state constraints and competing inputs was presented. This formulation was based on reachability and game theory, has the advantage of maintaining the continuity in the Hamiltonian of the system, and hence it has very good properties in terms of the numerical solution. The problem of reaching a desired set, in this case the TW, while avoiding conflict with other aircraft was formulated as a reach-avoid problem, and was computed numerically by using the existing tools.

In future work, we plan to use these reach-avoid bounds in order to perform conflict resolution by optimizing some cost criterion. Another issue would be to extend the proposed approach to formulate games in the case where the obstacle function is time and/or control dependant. Finally, we intend to validate our approach with fast time simulation studies using realistic aircraft and flight management system models, flight plans and wind uncertainty.

% or


\begin{appendix}
\section{}
\subsection{Proof of Proposition 2.}
\begin{proof}
\textbf{Part 1.} Following \cite{Mitchell_Bayen2} we first show that $\widetilde{RA}(\tau,R,A) \subseteq \{x \in \mathbb{R}^n ~|~ \widetilde{V}(x,\tau)\leq 0\}$. Consider $x \in \widetilde{RA}(\tau,R,A)$ and for the sake of contradiction assume that $\widetilde{V}(x,t)>0$.
Then there exists $\epsilon>0$ such that $\sup_{v(\cdot) \in  \mathcal{V}_{[t,T]}} \nonumber
  \max \{l(\tilde{\phi}(T,t,x,\tilde{\gamma}[v](\cdot),v(\cdot))) , \max_{\tau \in [t,T]} h(\tilde{\phi}(\tau,t,x,\tilde{\gamma}[v](\cdot),\\v(\cdot))) \} >2\epsilon > 0 $. This in turn implies that there exists $\hat{v}(\cdot) \in  \mathcal{V}_{[t,T]}$ such that
  either $l(\tilde{\phi}(T,t,x,\tilde{\gamma}[\hat{v}]\\(\cdot),\hat{v}(\cdot)))>\epsilon>0$ or there exists $\tau \in [t,T]$ such that $h(\tilde{\phi}(\tau,t,x,\tilde{\gamma}[\hat{v}](\cdot),\hat{v}(\cdot))) >\epsilon>0$.

Consider now the implications of $x \in \widetilde{RA}(\tau,R,A)$. Equation $(5)$ implies that there exists a $\gamma(\cdot) \in \Gamma_{[t,T]}$ such that for all $v(\cdot) \in \mathcal{V}_{[t,T]}$, and so also for $\hat{v}(\cdot)$, we can define $u(\cdot)=\gamma[\hat{v}](\cdot)$. Then, for this $u(\cdot)$ and $\hat{v}(\cdot)$ there exists $\tau_1 \in [t,T]$ such that $\phi(\tau_1,x,t,u(\cdot),\hat{v}(\cdot)) \in R$ and for all $\tau_2 \in [t,\tau_1] ~ \phi(\tau_2,t,x,u(\cdot),\hat{v}(\cdot)) \notin A$. Choose the freezing input signal as
\begin{equation*}
\bar{u}(s)= \left\{
\begin{array}{rl}
1 & \text{for } s \in [t,\tau_1]\\
0 & \text{for } s \in [\tau_1,T]
\end{array} \right.
\end{equation*}
If we combine $\bar{u}(\cdot)$ with $u(\cdot)$, we can get the input $\tilde{u}(\cdot)$ which will generate a trajectory
\begin{equation*}
\tilde{\phi}(\tau, x,t,\tilde{u}(\cdot),\hat{v}(\cdot))= \left\{
\begin{array}{rl}
\phi(\tau, x,t,u(\cdot),\hat{v}(\cdot)) & \text{for } \tau \in [t,\tau_1]\\
\phi(\tau_1, x,t,u(\cdot),\hat{v}(\cdot)) & \text{for } \tau \in [\tau_1,T]
\end{array} \right.
\end{equation*}
\textbf{Case 1.1:} Consider first the case where for all $\tilde{\gamma}(\cdot) \in \widetilde{\Gamma}_{[t,T]} ~l(\tilde{\phi}(T,t,x, \tilde{\gamma}[\hat{v}](\cdot),\hat{v}(\cdot))) >\epsilon>0$.
For $s=T$ we have that
\begin{equation*}
\tilde{\phi}(T, x,t,\tilde{u}(\cdot),\hat{v}(\cdot)) = \phi(\tau_1, x,t,u(\cdot),\hat{v}(\cdot)).
\end{equation*}
Since $x \in \widetilde{RA}(\tau,R,A)$, we showed before that $\phi(\tau_1,x,t,u(\cdot),\hat{v}(\cdot)) \in R$, i.e. $l(\phi(\tau_1,x,t,u(\cdot),\hat{v}(\cdot)\\))\leq 0$. So from we have that $l(\tilde{\phi}(T, x,t,\tilde{u}(\cdot),\hat{v}(\cdot))) \leq 0$. Since $u(\cdot)=\gamma[\hat{v}](\cdot)$ is already non-anticipative, and a non-anticipative strategy for $\bar{u}(\cdot)$ can be designed, $\tilde{u}(\cdot)$ will also be non-anticipative. Therefore, the previous statement establishes a contradiction. \\
\newline
\textbf{Case 1.2:} Consider now the case where for all $\tilde{\gamma}(\cdot) \in \widetilde{\Gamma}_{[t,T]}$ there exists $\tau \in [t,T]$ such that $h(\tilde{\phi}(\tau,t,x,\tilde{\gamma}[\hat{v}](\cdot),\hat{v}(\cdot))) >\epsilon>0$. Since we showed that for all $\tau \in [t,\tau_1], ~ \phi(\tau,t,x,u(\cdot),\hat{v}(\cdot)) \notin A$, we can conclude that for all $\tau \in [t,\tau_1]$
\begin{equation*}
h(\tilde{\phi}(\tau, x,t,\tilde{u}(\cdot),\hat{v}(\cdot))) \leq 0.
\end{equation*}
If $\tau \in [\tau_1,T]$, then we have that
\begin{equation*}
\tilde{\phi}(\tau, x,t,\tilde{u}(\cdot),\hat{v}(\cdot)) = \phi(\tau_1, x,t,u(\cdot),\hat{v}(\cdot)).
\end{equation*}
So $h(\tilde{\phi}(\tau, x,t,\tilde{u}(\cdot),\hat{v}(\cdot))) = h(\phi(\tau_1, x,t,u(\cdot),\hat{v}(\cdot))) \leq 0 $. Hence, for all $\tau \in [t,T]$ we have that $h(\tilde{\phi}(\tau, x,t,\tilde{u}(\cdot),\hat{v}(\cdot))) \leq 0$. Since in Case 1.1, $\tilde{u}(\cdot)$ was shown to be non-anticipative, we have a contradiction. \\
\newline
\textbf{Part 2.} Next, we show that $\{x \in \mathbb{R}^n ~|~ \widetilde{V}(x,\tau)\leq 0\} \subseteq \widetilde{RA}(\tau,R,A)$. Consider $(x,t)$ such that $\widetilde{V}(x,t) \leq 0$ and assume for the sake of contradiction that $x \notin \widetilde{RA}(\tau,R,A)$. Then for all for all $\gamma(\cdot) \in \Gamma_{[t,T]}$ there exists $\hat{v}(\cdot) \in \mathcal{V}_{[t,T]}$ such that for all $\tau_1 \in [t,T]$ either $\phi(\tau_1,t,x,\gamma(\cdot),v(\cdot)) \notin R$ or there exists $\tau_2 \in [t,\tau_1]$ such that $\phi(\tau_2,t,x,\gamma(\cdot),v(\cdot)) \in A$.

Following the analysis of \cite{Mitchell_Bayen2}, consider that the strategy $\gamma(\cdot) \in \Gamma_{[t,T]}$ is extracted from $\tilde{\gamma}(\cdot) \in \widetilde{\Gamma}_{[t,T]}$, and choose the $\hat{v}(\cdot)$ that corresponds to that strategy. In \cite{Mitchell_Bayen2}, was proven that the set of states visited by the augmented trajectory is a subset of the states visited by the original one. We therefore have that for all $\tau_1 \in [t,T]$
\begin{equation}\label{eq:15}
\phi(\tau_1, x,t,\gamma[\hat{v}](\cdot),\hat{v}(\cdot)) \notin R \Longrightarrow \tilde{\phi}(\tau_1, x,t,\tilde{\gamma}[\hat{v}](\cdot),\hat{v}(\cdot)) \notin R,
\end{equation}
and also $\forall \tau_2 \in [t,\tau_1]$
\begin{equation}\label{eq:16}
\phi(\tau_2, x,t,\gamma[\hat{v}](\cdot),\hat{v}(\cdot)) \in A \Longrightarrow \tilde{\phi}(\tau_2, x,t,\tilde{\gamma}[\hat{v}](\cdot),\hat{v}(\cdot)) \in A.
\end{equation}
By $(\ref{eq:15}),(\ref{eq:16})$ we conclude that there exists a $\delta > 0$ such that either for all $\tau_1 \in [t,T]$
\begin{equation}\label{eq:17}
l(\tilde{\phi}(\tau_1, x,t,\tilde{\gamma}[\hat{v}](\cdot),\hat{v}(\cdot))) > \delta > 0,
\end{equation}
or for some $\tau_2^* \in [t,\tau_1]$
\begin{equation}\label{eq:18}
h(\tilde{\phi}(\tau_2^*, x,t,\tilde{\gamma}[\hat{v}](\cdot),\hat{v}(\cdot))) > \delta > 0.
\end{equation}
Since $\widetilde{V}(x,t) \leq 0$, then for all $\epsilon>0$ there exists a non-anticipative strategy $\tilde{\gamma}(\cdot) \in \widetilde{\Gamma}_{[t,T]}$ such that $\sup_{v(\cdot) \in  \mathcal{V}_{[t,T]}},~ \max \{l(\tilde{\phi}(T,t,x,\tilde{\gamma}[v](\cdot),v(\cdot))) , \max_{\tau \in [t,T]} h(\tilde{\phi}(\tau,t,x,\tilde{\gamma}[v](\cdot),v(\cdot))) \} \leq \epsilon$. Hence for all $v(\cdot) \in  \mathcal{V}_{[t,T]}$, $l(\tilde{\phi}(T,t,x,\tilde{\gamma}[v](\cdot),v(\cdot)))\leq \epsilon$ and for all $\tau \in [t,T]$, $h(\tilde{\phi}(\tau,t,x,\tilde{\gamma}[v](\cdot),v(\cdot))) \leq \epsilon$.
For $v(\cdot)=\hat{v}(\cdot)$ the last argument implies that
\begin{equation*}
l(\tilde{\phi}(T, x,t,\tilde{\gamma}[\hat{v}](\cdot),\hat{v}(\cdot))) \leq \epsilon,
\end{equation*}
and there exists $\tau_2 \in [t,\tau_1]$
\begin{equation*}
h(\tilde{\phi}(\tau_2, x,t,\tilde{\gamma}[\hat{v}](\cdot),\hat{v}(\cdot))) \leq \epsilon.
\end{equation*}
If we choose $\epsilon = \frac{\delta}{2}$, the last statements contradict $(\ref{eq:17}),(\ref{eq:18})$ and complete the proof.
\end{proof}

\section{}
\subsection{Proof of Lemma 1.}
\begin{proof}
Following \cite{Evans_Souganidis} we can define
\begin{equation*} W(x,t) = \inf_{\gamma(\cdot)\in  \Gamma_{[t,t+\alpha]}} \sup_{v(\cdot)\in  \mathcal{V}_{[t,t+\alpha]}} \big[ \max \big \{ \max_{\tau \in [t,t+\alpha]} h(\phi(\tau,t,x,u(\cdot))), V(\phi(t+\alpha,t,x,u(\cdot)),t+\alpha) \big \} \big]. \end{equation*}
We will then show that for all $\epsilon > 0$, $V(x,t) \leq W(x,t) + 2\epsilon$ and $V(x,t) \geq W(x,t) - 3\epsilon$.
Then since $\epsilon > 0$ is arbitrary, $V(x,t)=W(x,t)$. \\
\newline
\textbf{Case 1: }$V(x,t) \leq W(x,t) + 2\epsilon$. Fix $\epsilon > 0$ and choose $\gamma_1(\cdot)\in  \Gamma_{[t,t+\alpha]}$ such that
\begin{align*} W(x,t) \geq \sup_{v_1(\cdot)\in  \mathcal{V}_{[t,t+\alpha]}} \big[ \max \big \{ &\max_{\tau \in [t,t+\alpha]} h(\phi(\tau,t,x,\gamma_1(\cdot),v_1(\cdot))), \\&V(\phi(t+\alpha,t,x,\gamma_1(\cdot),v_1(\cdot)),t+\alpha) \big \} \big] - \epsilon, \nonumber\end{align*}
Similarly, choose $\gamma_2(\cdot)\in  \Gamma_{[t+\alpha,T]}$ such that
\begin{align*}
V(\phi(t+\alpha,t,x,\gamma_1(\cdot),&v_1(\cdot)),t+\alpha) \geq \\&\sup_{v_2(\cdot)\in  \mathcal{V}_{[t+\alpha,T]}} \max \big \{ l(\phi(T,t+\alpha,\phi(t+\alpha,t,x,\gamma_1(\cdot),v_1(\cdot)),\gamma_2(\cdot),v_2(\cdot))), \\
&\max_{\tau \in [t+\alpha,T]} h(\phi(\tau,t+\alpha,\phi(t+\alpha,t,x,\gamma_1(\cdot),v_1(\cdot)),\gamma_2(\cdot),v_2(\cdot))) \big \}  - \epsilon.
\end{align*}
For any $v(\cdot)\in  \mathcal{V}_{[t,T]}$ we can define $v_1(\cdot)\in  \mathcal{V}_{[t,t+\alpha]}$ and $v_2(\cdot)\in  \mathcal{V}_{[t+\alpha,T]}$ such that $v_1(\tau)= v(\tau)$ for all $\tau \in [t,t+\alpha)$ and $v_2(\tau) = v(\tau)$ for all $\tau \in [t+\alpha,T]$.
Define also $\gamma(\cdot)\in  \Gamma_{[t,T]}$ by
\begin{equation*}
\gamma[v](\tau)= \left\{
\begin{array}{rl}
\gamma_1[v_1](\tau) & \text{if } \tau \in [t,t+\alpha)\\
\gamma_2[v_2](\tau) & \text{if } \tau \in [t+\alpha,T].
\end{array} \right.
\end{equation*}
It easy to see that $\gamma:~ \mathcal{V}_{[t,T]} \rightarrow \mathcal{U}_{[t,T]}$ is non-anticipative. By uniqueness, $\phi(\tau,t,x,\gamma(\cdot),v(\cdot))= \phi(\tau,t,x,\gamma_1(\cdot),v_1(\cdot))$ if $\tau \in [t,t+\alpha)$, and $\phi(\tau,t,x,\gamma(\cdot),v(\cdot))= \phi(\tau,t+\alpha,\phi(t+\alpha,t,x,\gamma_1(\cdot),v_1(\cdot)),\\\gamma_2(\cdot),v_2(\cdot))$ if $\tau \in [t+\alpha,T]$.\\
\newline
Hence,
\begin{align*}
W(x,t) &\geq  \sup_{v_1(\cdot)\in  \mathcal{V}_{[t,t+\alpha]}} \sup_{v_2(\cdot)\in  \mathcal{V}_{[t+\alpha,T]}} \max \big \{ \max_{\tau \in [t,t+\alpha]} h(\phi(\tau,t,x,\gamma_1(\cdot),v_1(\cdot))),
\\&l(\phi(T,t+\alpha,\phi(t+\alpha,t,x,\gamma_1(\cdot),v_1(\cdot)),\gamma_2(\cdot),v_2(\cdot))),\\
&\max_{\tau \in [t+\alpha,T]} h(\phi(\tau,t+\alpha,\phi(t+\alpha,t,x,\gamma_1(\cdot),v_1(\cdot)),\gamma_2(\cdot),v_2(\cdot))) \big \} - 2\epsilon\\
& \geq \sup_{v(\cdot)\in  \mathcal{V}_{[t,T]}}\max \big \{ l(\phi(T,t,x,\gamma(\cdot),v(\cdot))), \max_{\tau \in [t,T]} h(\phi(\tau,t,x,\gamma(\cdot),v(\cdot))) \big \} - 2\epsilon \\
& \geq V(x,t)-2\epsilon.
\end{align*}
Therefore, $V(x,t) \leq W(x,t) + 2\epsilon$. \\
\newline
\textbf{Case 2: }$V(x,t) \geq W(x,t) - 3\epsilon$. Fix $\epsilon > 0$ and choose now $\gamma(\cdot)\in  \Gamma_{[t,T]}$ such that
\begin{equation} \label{eq:7}
V(x,t) \geq \sup_{v(\cdot)\in  \mathcal{V}_{[t,T]}} \max \big \{ l(\phi(T,t,x,\gamma(\cdot),v(\cdot))), \max_{\tau \in [t,T]} h(\phi(\tau,t,x,\gamma(\cdot),v(\cdot))) \big \}-\epsilon.
\end{equation}
By the definition of $W(x,t)$
\begin{equation*} W(x,t) \leq \sup_{v(\cdot)\in  \mathcal{V}_{[t,t+\alpha]}} \big[ \max \big \{ \max_{\tau \in [t,t+\alpha]} h(\phi(\tau,t,x,\gamma(\cdot),v(\cdot))), V(\phi(t+\alpha,t,x,\gamma(\cdot),v(\cdot)),t+\alpha) \big \} \big].
\end{equation*}
Hence there exists a $v_1(\cdot) \in \mathcal{V}_{[t,t+\alpha]}$ such that
\begin{equation}\label{eq:8}
W(x,t) \leq \max \big \{ \max_{\tau \in [t,t+\alpha]} h(\phi(\tau,t,x,\gamma(\cdot),v_1(\cdot))), V(\phi(t+\alpha,t,x,\gamma(\cdot),v_1(\cdot)),t+\alpha) \big \} + \epsilon.
\end{equation}
Let $\hat{v}(\tau)= v_1(\tau)$ for all $\tau \in [t,t+\alpha)$ and $\hat{v}(\tau)= v'(\tau)$ for all $\tau \in [t+\alpha,T]$. Let also $\gamma' \in \Gamma_{[t+\alpha,T]}$ to be the restriction of the non-anticipative strategy $\gamma(\cdot)$ over $[t+\alpha,T]$. Then, for all $\tau \in [t+\alpha,T]$, we define $\gamma'[v'](\tau) = \gamma[\hat{v}](\tau)$.
Hence
\begin{align*}
V(\phi(t+\alpha,t,x,\gamma(\cdot),v_1(\cdot)),&t+\alpha) \leq \\&\sup_{v'(\cdot)\in  \mathcal{V}_{[t+\alpha,T]}} \max \big \{ l(\phi(T,t+\alpha,\phi(t+\alpha,t,x,\gamma(\cdot),v_1(\cdot)),\gamma'(\cdot),v'(\cdot))), \\
&\max_{\tau \in [t+\alpha,T]} h(\phi(\tau,t+\alpha,\phi(t+\alpha,t,x,\gamma(\cdot),v_1(\cdot)),\gamma'(\cdot),v'(\cdot))) \big \},
\end{align*}
and so there exists a $v_2(\cdot) \in \mathcal{V}_{[t+\alpha,T]}$ such that
\begin{align}\label{eq:9}
V(\phi(t+\alpha,t,x,\gamma(\cdot),v_1(\cdot)),&t+\alpha) \leq \max \big \{ l(\phi(T,t+\alpha,\phi(t+\alpha,t,x,\gamma(\cdot),v_1(\cdot)),\gamma'(\cdot),v_2(\cdot))), \\
&\max_{\tau \in [t+\alpha,T]} h(\phi(\tau,t+\alpha,\phi(t+\alpha,t,x,\gamma(\cdot),v_1(\cdot)),\gamma'(\cdot),v_2(\cdot))) \big \} + \epsilon. \nonumber
\end{align}
We can define
\begin{equation*}
v(\tau)= \left\{
\begin{array}{rl}
v_1(\tau) & \text{if } \tau \in [t,t+\alpha)\\
v_2(\tau) & \text{if } \tau \in [t+\alpha,T]
\end{array} \right.
\end{equation*}
Therefore, from (\ref{eq:8}) and (\ref{eq:9})
\begin{equation*}
W(x,t) \leq  \max \big \{ l(\phi(T,t,x,\gamma(\cdot),v(\cdot))), \max_{\tau \in [t,T]} h(\phi(\tau,t,x,\gamma(\cdot),v(\cdot))) \big \} + 2\epsilon, \nonumber
\end{equation*}
which together with (\ref{eq:7}) implies $V(x,t) \geq  W(x,t) - 3\epsilon$.
\end{proof}
\subsection{Proof of Lemma 2.}
\begin{proof}
Since $l$, and $h$ are bounded, $V$ is also bounded. For the second part fix $x, \hat{x} \in \mathbb{R}^n$ and $t \in [0,T]$. Let $\epsilon > 0$ and choose $\hat{\gamma}(\cdot) \in \Gamma_{[t,T]}$ such that
\begin{equation*}
V(\hat{x},t) \geq \sup_{v(\cdot) \in \mathcal{V}_{[t,T]}} \max_{\tau \in [t,T]} \max \{l(\phi(T,t,\hat{x},\hat{\gamma}(\cdot),v(\cdot))) , h(\phi(\tau,t,\hat{x},\hat{\gamma}(\cdot),v(\cdot))) \} - \epsilon.
\end{equation*}
By definition
\begin{equation*}
V(x,t) \leq \sup_{v(\cdot) \in \mathcal{V}_{[t,T]}} \max_{\tau \in [t,T]} \max \{l(\phi(T,t,x,\hat{\gamma}(\cdot),v(\cdot))) , h(\phi(\tau,t,x,\hat{\gamma}(\cdot),v(\cdot))) \}.
\end{equation*}
We can choose $\hat{v}(\cdot) \in \mathcal{V}_{[t,T]}$ such that
\begin{equation*}
V(x,t) \leq \max_{\tau \in [t,T]} \max \{l(\phi(T,t,x,\hat{\gamma}(\cdot),\hat{v}(\cdot))) , h(\phi(\tau,t,x,\hat{\gamma}(\cdot),\hat{v}(\cdot))) \} + \epsilon,
\end{equation*}
and hence
\begin{align*}
V(x,t)-V(\hat{x},t) &\leq \max_{\tau \in [t,T]} \max \{l(\phi(T,t,x,\hat{\gamma}(\cdot),\hat{v}(\cdot))) , h(\phi(\tau,t,x,\hat{\gamma}(\cdot),\hat{v}(\cdot))) \} \\ &- \max_{\tau \in [t,T]} \max \{l(\phi(T,t,\hat{x},\hat{\gamma}(\cdot),\hat{v}(\cdot))) , h(\phi(\tau,t,\hat{x},\hat{\gamma}(\cdot),\hat{v}(\cdot))) \} + 2\epsilon.
\end{align*}
For all $\tau \in [t,T]$:
\begin{align*}
|\phi(\tau,t,x,\hat{\gamma}(\cdot),\hat{v}(\cdot)) &- \phi(T,t,\hat{x},\hat{\gamma}(\cdot),\hat{v}(\cdot))| = \\& |(x-\hat{x})+\int_t^T [f(\phi(s,t,x,\hat{\gamma}(\cdot),\hat{v}(\cdot))) - f(\phi(s,t,\hat{x},\hat{\gamma}(\cdot),\hat{v}(\cdot)))]ds| \\
& \leq |x-\hat{x}|+\int_t^T |f(\phi(s,t,x,\hat{\gamma}(\cdot),\hat{v}(\cdot))) - f(\phi(s,t,\hat{x},\hat{\gamma}(\cdot),\hat{v}(\cdot)))|ds \\
& \leq |x-\hat{x}|+ C_f\int_t^T |\phi(s,t,x,\hat{\gamma}(\cdot),\hat{v}(\cdot)) - \phi(s,t,\hat{x},\hat{\gamma}(\cdot),\hat{v}(\cdot))|ds,
\end{align*}
where $C_f$ is the Lipschitz constant of $f$. By the Gronwall-Bellman Lemma \cite{Sastry}, there exists a constant $C_x > 0$ such that for all $\tau \in [t,T]$
\begin{equation*}
|\phi(\tau,t,x,\hat{\gamma}(\cdot),\hat{v}(\cdot)) - \phi(T,t,\hat{x},\hat{\gamma}(\cdot),\hat{v}(\cdot))| \leq C_x|x-\hat{x}|.
\end{equation*}
Let $\tau_0 \in [t,T]$ be such that
\begin{equation*}
h(\phi(\tau_0,t,x,\hat{\gamma}(\cdot),\hat{v}(\cdot))) = \max_{\tau \in [t,T]} h(\phi(\tau,t,x,\hat{\gamma}(\cdot),\hat{v}(\cdot))).
\end{equation*}
Then
\begin{align*}
V(x,t)-V(\hat{x},t) & \leq \max \{l(\phi(T,t,x,\hat{\gamma}(\cdot),\hat{v}(\cdot))) , h(\phi(\tau_0,t,x,\hat{\gamma}(\cdot),\hat{v}(\cdot))) \}\\ & - \max \{l(\phi(T,t,\hat{x},\hat{\gamma}(\cdot),\hat{v}(\cdot))) , h(\phi(\tau_0,t,\hat{x},\hat{\gamma}(\cdot),\hat{v}(\cdot))) \} + 2\epsilon.
\end{align*}
\textbf{Case 1.} $~l(\phi(T,t,x,\hat{\gamma}(\cdot),\hat{v}(\cdot))) \geq h(\phi(\tau_0,t,x,\hat{\gamma}(\cdot),\hat{v}(\cdot)))$
\begin{align*}
V(x,t)-V(\hat{x},t) &\leq l(\phi(T,t,x,\hat{\gamma}(\cdot),\hat{v}(\cdot))) \\&- \max \{l(\phi(T,t,\hat{x},\hat{\gamma}(\cdot),\hat{v}(\cdot))) , h(\phi(\tau_0,t,\hat{x},\hat{\gamma}(\cdot),\hat{v}(\cdot))) \} + 2\epsilon \\
& \leq l(\phi(T,t,x,\hat{\gamma}(\cdot),\hat{v}(\cdot))) - l(\phi(T,t,\hat{x},\hat{\gamma}(\cdot),\hat{v}(\cdot))) + 2\epsilon \\
& \leq C_l C_x |x-\hat{x}| + 2\epsilon.
\end{align*}
\textbf{Case 2.} $~l(\phi(T,t,x,\hat{\gamma}(\cdot),\hat{v}(\cdot))) < h(\phi(\tau_0,t,x,\hat{\gamma}(\cdot),\hat{v}(\cdot)))$
\begin{align*}
V(x,t)-V(\hat{x},t) &\leq h(\phi(\tau_0,t,x,\hat{\gamma}(\cdot),\hat{v}(\cdot))) \\&- \max \{l(\phi(T,t,\hat{x},\hat{\gamma}(\cdot),\hat{v}(\cdot))) , h(\phi(\tau_0,t,\hat{x},\hat{\gamma}(\cdot),\hat{v}(\cdot))) \} + 2\epsilon \\
& \leq h(\phi(\tau_0,t,x,\hat{\gamma}(\cdot),\hat{v}(\cdot))) - h(\phi(\tau_0,t,\hat{x},\hat{\gamma}(\cdot),\hat{v}(\cdot))) + 2\epsilon \\
& \leq C_h C_x |x-\hat{x}| + 2\epsilon.
\end{align*}
So in any case $ V(x,t)-V(\hat{x},t)\leq \max \{C_l C_h\} C_x |x-\hat{x}| + 2\epsilon$. The same argument with the roles of $x$, $\hat{x}$ reversed establishes that $ V(\hat{x},t)-V(x,t)\leq \max \{C_l C_h\} C_x |x-\hat{x}| + 2\epsilon$. Since $\epsilon$ is arbitrary,
\begin{equation*}
|V(x,t)-V(\hat{x},t)| \leq \max \{C_l C_h\} C_x |x-\hat{x}|.
\end{equation*}
Finally consider $x \in \mathbb{R}^n$ and $t, \hat{t} \in [0,T]$. Without loss of generality assume that $t<\hat{t}$.
Let $\epsilon > 0$ and choose $\gamma(\cdot) \in \Gamma_{[t,T]}$ such that
\begin{align*}
V(x,t) &\geq \sup_{v(\cdot) \in \mathcal{V}_{[t,T]}} \max_{\tau \in [t,T]} \max \{l(\phi(T,t,x,\gamma(\cdot),v(\cdot))), h(\phi(\tau,t,x,\gamma(\cdot),v(\cdot))) \} - \epsilon \\
&\geq \max_{\tau \in [t,T]} \max \{l(\phi(T,t,x,\gamma(\cdot),v(\cdot))), h(\phi(\tau,t,x,\gamma(\cdot),v(\cdot))) \} - \epsilon
\end{align*}
By definition,
\begin{equation*}
V(x,\hat{t}) \leq \sup_{v(\cdot) \in \mathcal{V}_{[\hat{t},T]}} \max_{\tau \in [\hat{t},T]} \max \{l(\phi(T,\hat{t},x,\hat{\gamma}(\cdot),v(\cdot))) , h(\phi(\tau,\hat{t},x,\hat{\gamma}(\cdot),v(\cdot))) \}.
\end{equation*}
So we can choose $\hat{v}(\cdot) \in \mathcal{V}_{[\hat{t},T]}$ such that
\begin{equation*}
V(x,\hat{t}) \leq \max_{\tau \in [\hat{t},T]} \max \{l(\phi(T,\hat{t},x,\hat{\gamma}(\cdot),\hat{v}(\cdot))) , h(\phi(\tau,\hat{t},x,\hat{\gamma}(\cdot),\hat{v}(\cdot))) \} + \epsilon,
\end{equation*}
where $\hat{\gamma} \in \Gamma_{[\hat{t},T]}$ is the restriction of $\gamma(\cdot)$ over $[\hat{t},T]$. Then, for all $\tau \in [\hat{t},T]$, we define $\hat{\gamma}[\hat{v}](\tau) = \gamma[v](\tau)$, and $\hat{v}(\tau) = v(\tau+t-\hat{t})$. By uniqueness, for all $\tau \in [\hat{t},T]$ we have that $\phi(\tau,\hat{t},x,\hat{\gamma}(\cdot),\hat{v}(\cdot))=\phi(\tau+t-\hat{t},t,x,\gamma(\cdot),v(\cdot))$.
\begin{align*}
V(x,t)-V(x,\hat{t}) &\geq \max_{\tau \in [t,T]} \max \{l(\phi(T,t,x,\gamma(\cdot),v(\cdot))) , h(\phi(\tau,t,x,\gamma(\cdot),v(\cdot))) \}\\ & - \max_{\tau \in [\hat{t},T]} \max \{l(\phi(T,\hat{t},x,\hat{\gamma}(\cdot),\hat{v}(\cdot))) , h(\phi(\tau,\hat{t},x,\hat{\gamma}(\cdot),\hat{v}(\cdot))) \} - 2\epsilon.
\end{align*}
\textbf{Case 1.} $~l(\phi(T,\hat{t},x,\hat{\gamma}(\cdot),\hat{v}(\cdot))) \geq \max_{\tau \in [\hat{t},T]} h(\phi(\tau,\hat{t},x,\hat{\gamma}(\cdot),\hat{v}(\cdot)))$
\begin{align*}
V(x,t)-V(x,\hat{t}) &\geq \max_{\tau \in [t,T]} \max \{l(\phi(T,t,x,\gamma(\cdot),v(\cdot))) , h(\phi(\tau,t,x,\gamma(\cdot),v(\cdot))) \} \\&-l(\phi(T,\hat{t},x,\hat{\gamma}(\cdot),\hat{v}(\cdot)))  - 2\epsilon \\
& \geq l(\phi(T,t,x,\gamma(\cdot),v(\cdot))) - l(\phi(T,\hat{t},x,\hat{\gamma}(\cdot),\hat{v}(\cdot))) - 2\epsilon \\
& = l(\phi(T,t,x,\gamma(\cdot),v(\cdot))) - l(\phi(T+t-\hat{t},t,x,\gamma(\cdot),v(\cdot))) - 2\epsilon \\
& \geq -C_l C_f |T-T-t+\hat{t}|-2\epsilon \\
& = -C_l C_f |\hat{t}-t| - 2\epsilon,
\end{align*}
where $C_l$ is the Lipschitz constant of $l$. \\
\newline
\textbf{Case 2.} $~l(\phi(T,\hat{t},x,\hat{\gamma}(\cdot),\hat{v}(\cdot))) < \max_{\tau \in [\hat{t},T]} h(\phi(\tau,\hat{t},x,\hat{\gamma}(\cdot),\hat{v}(\cdot)))$
\begin{align*}
V(x,t)-V(x,\hat{t}) &\geq \max_{\tau \in [t,T]} \max \{l(\phi(T,t,x,\gamma(\cdot),v(\cdot))) , h(\phi(\tau,t,x,\gamma(\cdot),v(\cdot))) \} \\&-\max_{\tau \in [\hat{t},T]} h(\phi(\tau,\hat{t},x,\hat{\gamma}(\cdot),\hat{v}(\cdot)))  - 2\epsilon \\
& \geq \max_{\tau \in [t,T]} h(\phi(\tau,t,x,\gamma(\cdot),v(\cdot))) -\max_{\tau \in [\hat{t},T]} h(\phi(\tau,\hat{t},x,\hat{\gamma}(\cdot),\hat{v}(\cdot)))  - 2\epsilon.
\end{align*}
Let $\tau_0 \in [\hat{t},T]$ be such that
\begin{equation*}
h(\phi(\tau_0,\hat{t},x,\hat{\gamma}(\cdot),\hat{v}(\cdot))) = \max_{\tau \in [\hat{t},T]} h(\phi(\tau,\hat{t},x,\hat{\gamma}(\cdot),\hat{v}(\cdot))).
\end{equation*}
Then
\begin{align*}
V(x,t)&-V(x,\hat{t}) \geq \max_{\tau \in [t,T]} h(\phi(\tau,t,x,\gamma(\cdot),v(\cdot))) -h(\phi(\tau_0,\hat{t},x,\hat{\gamma}(\cdot),\hat{v}(\cdot)))  - 2\epsilon \\
& \geq h(\phi(\tau_0,t,x,\gamma(\cdot),v(\cdot))) -h(\phi(\tau_0,\hat{t},x,\hat{\gamma}(\cdot),\hat{v}(\cdot)))  - 2\epsilon \\
& = h(\phi(\tau_0,t,x,\gamma(\cdot),v(\cdot))) -h(\phi(\tau_0+t-\hat{t},t,x,\gamma(\cdot),v(\cdot)))  - 2\epsilon \\
& \geq -C_h C_f |\tau_0-\tau_0-t+\hat{t}|-2\epsilon \\
& = -C_h C_f |\hat{t}-t|-2\epsilon,
\end{align*}
where $C_h$ is the Lipschitz constant of $h$. In any case we have that
\begin{equation*}
V(x,t)-V(x,\hat{t}) \geq -\max\{C_l,C_h\} C_f |\hat{t}-t|-2\epsilon.
\end{equation*}
A symmetric argument shows that $V(x,t)-V(x,\hat{t}) \leq \max\{C_l,C_h\} C_f |\hat{t}-t|+2\epsilon$, and since $\epsilon$ is arbitrary this concludes the proof.
\end{proof}
%\section{}

\section{}
The following algorithm summarizes the steps of the Reach-Avoid computation described in Section IV. For simplicity, we have assumed that the TW do not overlap.
\begin{algorithm}[h!]
\label{alg:SCPF}

\caption{\textbf{Reach-Avoid computation}}

\begin{algorithmic}[1]

\vspace{0.5cm}

\STATE \textbf{Initialization}. Set:\\
$T = \max_j \overline{t_j},$ for $j=1,...,N$,\\
$t_0 = \min_j \underline{t_j},$ for $j=1,...,N$,\\
{\bf  $\widetilde{V}(x_j,\overline{t}_j) = l(x_j)$ },\\
{\bf $\widetilde{R}_j = \{x_j ~|~ \widetilde{V}(x_j,\overline{t}_j) \leq 0 \}$ }.

\FOR{$t=T$ until $t_0$}

\STATE \textbf{while $j$ is in the sector}

\STATE ~~~~\textbf{if $\underline{t}_j \leq t \leq \overline{t}_j$} \\
~~~~~~~~Solve $\frac{\partial \widetilde{V}}{\partial t}(x_j,t)+ \min \{0,\sup_{v \in \textit{V}} \inf_{u_j \in \textit{U}_j} \frac{\partial \widetilde{V}}{\partial x_j}(x_j,t) f_j(x_j,u_j,v)\}=0$.

\STATE ~~~~~~~~\textbf{for all $i \neq j$ in the sector}\\
~~~~~~~~~~~~$C_{ji} = \{x_j ~|~ i,j \text{ are in conflict}\}$.

 \STATE ~~~~~~~~~~~~\textbf{for all $i,j \in C_{ji}$}\\
~~~~~~~~~~~~~~~~Define $h_{ji}(x_j)$ such that $A_{ji}\supseteq C_{ji}$ is a box.\\
~~~~~~~~~~~~~~~~Let $A_{ji} = \{x_j ~|~ h_{ji}(x_j)>0\}$, and $A_{j} = \bigcup_{i \neq j}A_{ji}$, \\
~~~~~~~~~~~~~~~~$h_{j} = \max_{i \neq j} h_{ji}(x_j)$,\\
~~~~~~~~~~~~~~~~$\widetilde{V}(x_j,t) = \max(h_{j},\widetilde{V}(x_j,t))$.\\
~~~~~~~~~~~~\textbf{else}\\
~~~~~~~~~~~~~~~~$A_{j} = \emptyset$. \\
\STATE~~~~~~~~~~~~\textbf{end for}\\
\STATE~~~~~~~~\textbf{end for}\\
\STATE~~~~~~~~\textbf{$V(x_j,\underline{t}_j) = \widetilde{V}(x_j,t)$}.\\
\STATE ~~~~\textbf{else if $t \leq \underline{t}_j$} \\
~~~~~~~~Solve $\frac{\partial V}{\partial t}(x_j,t)+ \sup_{v \in \textit{V}} \inf_{u_j \in \textit{U}_j} \frac{\partial V}{\partial x_j}(x_j,t) f_j(x_j,u_j,v)=0$.\\  
~~~~~~~~Repeat steps $5-8$ with $V$ instead of $\widetilde{V}$.\\
\STATE ~~~~\textbf{end if} \\
\STATE \textbf{end while}
\ENDFOR
\vspace{0.5cm}

\end{algorithmic}
\end{algorithm}
\newpage
\end{appendix}
\section*{Acknowledgment}
Research was supported by the European Commission under the project CATS, FP6-TREN-036889.

%\bibliographystyle{IEEEtran}
%\bibliography{refReachAvoid}


\end{document}